\documentclass[11pt,a4paper]{article}
\usepackage{amstex}
\usepackage{amssymb}
\newcommand{\gM}{{\mathfrak M}}
\newcommand{\gN}{{\mathfrak N}}
\newcommand{\Def}{{\it Def}}
\newcommand{\rpfeil}[2]{\stackrel{#2}{\verylongarrow{#1mm}}}
\newcommand{\dual}{\makebox[0mm]{}^{{\scriptstyle\vee}}}
\newtheorem{theorem}{Theorem}[section]
\newtheorem{lemma}[theorem]{Lemma}
\newtheorem{proposition}[theorem]{Proposition}
\newtheorem{definition}[theorem]{Definition}
\newtheorem{corollary}[theorem]{Corollary}
\newtheorem{exmp}[theorem]{Example}
\newtheorem{exmps}[theorem]{Examples}
\newtheorem{rem}[theorem]{Remark}
\newenvironment{example}{\begin{exmp}\rm}{\end{exmp}}

\newenvironment{remark}{\begin{rem}\rm}{\end{rem}\rm}
\newcommand{\prf}{{\em Proof}. }
\newcommand{\qed}{\hspace*{\fill}$\Box$}

\newcommand{\beeq}[1]{\begin{eqnarray}\label{#1}}
\newcommand{\eneq}{\end{eqnarray}}

\newcommand{\ka}{{\cal A}}

\newcommand{\kc}{{\cal C}}
\newcommand{\kd}{{\cal D}}
\newcommand{\kh}{{\cal H}}
\newcommand{\kk}{{\cal K}}
\newcommand{\kl}{{\cal L}}
\newcommand{\kn}{{\cal N}}
\newcommand{\ko}{{\cal O}}
\newcommand{\kp}{{\cal P}}
\newcommand{\kt}{{\cal T}}
\newcommand{\kx}{{\cal X}}

\newcommand{\IC}{{\mathbb C}} 
\newcommand{\IH}{{\mathbb H}}

\newcommand{\IP}{{\mathbb P}} 
\newcommand{\IQ}{{\mathbb Q}} 
\newcommand{\IR}{{\mathbb R}}
\newcommand{\IZ}{{\mathbb Z}}

\newcommand{\Hilb}{{\rm Hilb}}
\newcommand{\Pic}{{\rm Pic}}
\newcommand{\Hom}{{\rm Hom}}

\newcommand{\id}{{\rm id}}
\newcommand{\im}{{\rm im}}

\newcommand{\codim}{{\rm codim}}

\newcommand{\verylongarrow}[1]{\hbox to #1{\rightarrowfill}}

\begin{document}
{\parindent0mm{\Large\bf Compact Hyperk\"ahler Manifolds: Basic Results}}

\bigskip
\bigskip

{\parindent0mm{\large\bf Daniel Huybrechts}}\\
{\small Universit\"at--GH Essen, Fachbereich 6 Mathematik, 45141
Essen, Germany\\ daniel.huybrechts@@uni-essen.de}

\bigskip

\bigskip

\bigskip

\bigskip

Compact hyperk\"ahler manifolds, or irreducible symplectic manifolds
as they will be frequently called in these notes, are higher-dimensional
analogues of K3 surfaces. That they indeed share many of the
well-known properties of K3 surfaces this paper intends to show.

To any K3 surface $X$ one naturally associates its period consisting
of the weight-two Hodge structure on $H^2(X,\IZ)$ together
with the intersection pairing.
The period describes and determines a K3 surface
in the following sense:

\bigskip
\begin{itemize}
\item Firstly, the period encodes the algebraic (and to some
extent also the differential)
geometry of the K3 surface. E.g.\ the ample (K\"ahler) cone
can be described purely in terms of the period.
\item Secondly, any point in the period domain is realized as the period
of a K3 surface.
\item Thirdly, and this is certainly the culmination of the theory,
two K3 surfaces with the same period are isomorphic (Global Torelli
Theorem).
\end{itemize}

\bigskip

The second cohomology $H^2(X,\IZ)$ of a compact irreducible
hyperk\"ahler (or irreducible symplectic) manifold also carries a natural
weight-two Hodge structure. Moreover, due to a result of Beauville $H^2(X,\IZ)$
can be endowed with a natural non-degenerate quadratic form $q_X$ generalizing
the intersection pairing of a K3 surface.
For higher-dimensional varieties, in contrast to varieties of dimension
two, the knowledge of the second cohomology usually reveals only a small
portion of the full geometry. However, for irreducible
symplectic manifolds the situation is quite different and it is
not completely unreasonable to expect $(H^2(X,\IZ),q_X)$
to encode all `essential' information about the geometry of $X$.
In this paper we present several results underpinning this
expectation.
But to temper the enthusiasm a little at this point we mention that
Debarre has constructed examples showing that the Global Torelli
Theorem as formulated for K3 surfaces fails in higher dimensions.

\bigskip
The first basic result showing that $(H^2(X,\IZ),q_X)$
controls the geometry of $X$ at least to a certain extent
is the following projectivity criterion (Theorem \ref{proj})
which is also of importance for many of the other results in the
later chapters.

\bigskip
\begin{itemize}
\item An irreducible symplectic manifold $X$ is projective
if and only if there exists a line bundle $L$ on $X$ with
positive square $q_X(c_1(L))>0$.
\end{itemize}

\bigskip
Another instance of the significance of the second
cohomology for the geometry of $X$ is a description of the ample cone.
Due to a missing concept of $(-2)$-classes in higher dimensions the result is
less explicit than the known one for K3 surfaces. It can roughly be
formulated as

\bigskip
\begin{itemize}
\item A class in $H^2(X,\IZ)$ is ample if and only if it
has positive square with respect to $q_X$  and
cannot be separated from the K\"ahler cone by any (integral) hyperplane.
\end{itemize}

\bigskip
For the precise statement see Theorem \ref{ampleconethm}. A similar description
is obtained for the K\"ahler cone (Theorem
\ref{kaehlerconethm}). Again, the result does not reach the
explicitness of the result for K3 surfaces, but it suffices
to describe the K\"ahler cone of certain hyperk\"ahler manifolds
completely. One of the consequences is the following result
(Corollary \ref{kaehlerconecor}):

\bigskip
\begin{itemize}
\item For a general irreducible symplectic manifold $X$ the
K\"ahler cone equals the positive cone, i.e.\ the K\"ahler cone is a connected
component of $\{\alpha\in H^{1,1}(X,\IR)|q_X(\alpha)>0\}$.
\end{itemize}

\bigskip
An irreducible symplectic manifold is general if the corresponding
point in the moduli space is contained in the complement
of countably many nowhere dense closed subspaces (see \ref{deformationallg}).
In fact, the description of the K\"ahler cone is based on 
the higher dimensional analogue of the transitivity of the Weyl-action
for K3 surfaces.
The latter says that any general $(1,1)$-class with positive selfintersection
on a K3 surface can be mapped to a K\"ahler class by reflecting it in a
finite number of hyperplanes orthogonal to $(-2)$-curves.
Proposition \ref{Weylprop} and Theorem \ref{Delignesuggestion} are in this
spirit, although in higher dimensions there does
not exist a Weyl-group and the resulting class is only K\"ahler
on an irreducible symplectic manifold birational to the given one.

The information about the
K\"ahler cone allows us to attack the question of which point in the
period domain really occurs as the period of an irreducible
symplectic manifold:

\bigskip
\begin{itemize}
\item The period map surjects any non-empty connected component
of the moduli space of marked irreducible symplectic manifolds onto the
period domain $\{x|q_\Gamma(x)=0,~q_\Gamma(x+\bar x)>0\}\subset \IP(\Gamma_\IC)$.
\end{itemize}

\bigskip

This is Theorem \ref{surjper} (cf.\ Sect.\ \ref{prel} for the notation).
The last result we wish to mention explicitly in the
introduction deals with birational irreducible
symplectic manifolds. It has been known for some time that
seemingly different examples of irreducible symplectic manifolds
are often related to each other
either by deformation or by birational correspondences.
In \cite{Huybrechts} we could show that certain moduli spaces of stable
rank two sheaves on a K3 surface
are deformation equivalent (and hence diffeomorphic) to the Hilbert
scheme of points on this K3 surface. In fact, we showed
more generally that under a certain assumption on the codimension
of the exceptional sets any two birational irreducible
symplectic manifolds are diffeomorphic. Using the techniques developed
in this paper we are now able to prove this result
in full generality. As a consequence of the main Theorem \ref{birat}
we mention:

\bigskip
\begin{itemize}
\item Two birational projective irreducible symplectic manifolds
are deformation equivalent and, hence, diffeomorphic.
\end{itemize}

\bigskip
This replaces the fact that birational K3 surfaces are isomorphic.
The result also compares nicely with a recent theorem of Batyrev and Kontsevich
saying in particular that birational projective manifolds with trivial
canonical bundle have the same Betti numbers (see the discussion in
Sect.\ \ref{biratman}).

\bigskip
The plan of the paper is as follows:
In Sect.\ \ref{prel} we recall the basic concepts and
formulate some of the known results in a form ready for later use.
Some of the statements I could not find
in the literature, though they are probably known to the experts.
Sect.\ 1 has grown out of size a little, but I hope it makes the rest
of the paper more easily accessible.
Sect.\ \ref{examples} gives a list of all examples (known to
me) of irreducible symplectic manifolds.
This is important for a result on the classification
of these examples in Sect.\ \ref{biratman}.
Sect.\ \ref{projectivity} deals
with the above mentioned projectivity criterion.
In Sect.\ \ref{biratman} we show that the Main Lemma of Burns and
Rapoport \cite{BurnsRapoport} can be adapted easily to
the higher-dimensional situation.
Then we generalize the main result of \cite{Huybrechts}
and show that two birational projective irreducible symplectic
manifolds correspond to non-separated points on the moduli space. This
immediately implies the fact about the deformation equivalence
as stated above. This result is then applied to the list of examples
in Sect.\ \ref{examples}.
In Sect.\ \ref{weylaction} we present a result which can be seen
as a replacement of the transitivity of the action of the Weyl-group
on the set of chambers of a K3 surface. It has interesting
applications to the K\"ahler cone (Sect.\ \ref{kaehlercone})
and leads eventually to 
the surjectivity of the period map in Sect.\ \ref{periodmap}.
Compared to the standard proofs for K3 surfaces the order of arguments has
got an interesting twist.
The ampleness of line bundles on irreducible symplectic manifolds
is discussed in Sect.\ \ref{amplecone}. The penultimate
Sect.\ \ref{auto} is devoted to the
automorphism group of irreducible symplectic manifolds. We collect
known results and formulate some open problems.
The results of Sect.\ \ref{projectivity}-\ref{periodmap} clearly indicate
that irreducible symplectic manifolds behave in many respects like K3 surfaces.
However, the two main problems in the theory remain wide open and
seem to require completely different techniques (this also explains
the `Basic Results' in the title): {\it i)}
Is there any kind of Global Torelli
Theorem for higher-dimensional irreducible symplectic
manifolds? and {\it ii)} What are the possible deformation (diffeomorphism,
homeomorphism) types of irreducible symplectic manifolds?
We discuss in Sect.\ \ref{remarks} some of the aspects of
these two problems and
their relation to questions on the moduli space of marked manifolds.

\bigskip
{\bf Acknowledgements:} My main source of inspiration were Beauville's
presentation of the theory of K3 surfaces in \cite{Periodes} and Fujiki's
papers on symplectic manifolds \cite{Fujiki1,Fujiki2}.
The bulk of this paper was written while I was enjoying the hospitality
of the IHES (Bures-sur-Yvette) and the Humboldt-Universit\"at (Berlin).
I am also grateful for financial support provided by the Max-Planck-Institut
(Bonn), the IHES, and the Universit\"at-GH Essen.
Several people have read preliminary versions of this article and
made valuable comments. In particular, I express my gratitude
to L.\ Bonavero, P.\ Deligne,\ H.\ Esnault, and E.\ Viehweg.
P.\ Deligne suggested an improvement of Proposition \ref{Weylprop}
which resulted in Theorem \ref{Delignesuggestion}.
I wish to thank H.\ Esnault and E.\ Viehweg also for their warm welcome
to Essen and their encouragement in the final stage of this work.

\section{Preliminaries}\label{prel}

{\bf \refstepcounter{theorem}\label{IrredSympl} \thetheorem} ---
A complex manifold $X$ is called {\it irreducible symplectic} if
{\it i)} $X$ is compact and K\"ahler,
{\it ii)} $X$ is simply connected, and {\it iii)}
$H^0(X,\Omega^2_X)$ is spanned by an everywhere
non-degenerate two-form $\sigma$.

Any holomorphic two-form $\sigma$ induces a homomorphism $\kt_X\to \Omega_X$,
which we also denote by $\sigma$. The two-form is everywhere non-degenerate
if and only if $\sigma:\kt_X\to\Omega_X$ is bijective. Note that {\it iii)}
implies $h^{2,0}(X)=h^{0,2}(X)=1$ and $K_X\cong \ko_X$. In particular,
$c_1(X)=0$. Any irreducible symplectic manifold $X$ has even complex dimension
which we will fix to be $2n$. Irreducible symplectic
manifolds occupy a distinguished place in the list of higher dimensional
K\"ahler manifolds. Together with Calabi-Yau manifolds they are the only
irreducible simply connected K\"ahler manifolds with $c_1(X)=0$ (cf.\
\cite{Beauville1}).

\bigskip
{\bf \refstepcounter{theorem}\label{HypK} \thetheorem} ---
A compact connected $4n$-dimensional Riemannian manifold $(M,g)$ is called
{\it hyperk\"ahler} ({\it irreducible hyperk\"ahler}) if its holonomy
is contained in (equals) ${\rm Sp}(n)$.

If $(M,g)$ is hyperk\"ahler, then the quaternions $\IH$ act as
parallel endomorphisms on the tangent bundle of $M$. This is a consequence
of the holonomy principle: Every tensor at a point in $M$ that is invariant
under the holonomy action can be extended to a parallel tensor over $M$.
In particular, any $\lambda\in\IH$ with $\lambda^2=-1$ gives rise to an almost
complex structure on $M$. As it turns out, these almost complex structures
are all integrable \cite{Salamon2}. After having fixed a standard basis
$I$, $J$, and $K:=IJ$ of $\IH$ any
$\lambda\in\IH$ with $\lambda^2=-1$ can be written as $\lambda=aI+bJ+cK$
with $a^2+b^2+c^2=1$. Note that the metric
$g$ is K\"ahler with respect to 
every such $\lambda\in S^2$. The corresponding K\"ahler form is denoted by
$\omega_\lambda:=g(\lambda\,.\,,\,.\,)$.

Thus, a hyperk\"ahler metric $g$ on a manifold $M$ defines a family
of complex K\"ahler manifolds $(M,\lambda,\omega_\lambda)$,
where $\lambda\in S^2\cong\IP^1$.

\bigskip
{\bf \refstepcounter{theorem}\label{Comparison} \thetheorem} ---
We briefly sketch the relation between irreducible symplectic and
irreducible hyperk\"ahler manifolds. For details we refer to
\cite{Beauville1}. 

If $X$ is irreducible symplectic and $\alpha\in H^2(X,\IR)$ is a K\"ahler class
on $X$, then there exists a unique Ricci-flat K\"ahler metric $g$ with
K\"ahler class $\alpha$. This follows from Yau's solution of the
Calabi-conjecture. Then $g$ is an irreducible hyperk\"ahler metric
on the underlying real manifold $M$. Moreover, for one of the
complex structures, say $I$, one has $X=(M,I)$.

Conversely, if $(M,g)$ is hyperk\"ahler and $I,J,K$ are complex structures
as above, then $\sigma:=\omega_J+i\omega_K$ is a holomorphic everywhere
non-degenerate two-form
on $X=(M,I)$. If $M$ is compact and $g$ is irreducible hyperk\"ahler, then
$M$ is simply connected and $H^0((M,I),\Omega_{(M,I)}^2)=\sigma\cdot\IC$,
i.e.\ $X$ is irreducible symplectic.

Thus, irreducible symplectic manifolds with a fixed K\"ahler class and
compact irreducible hyperk\"ahler manifolds are the holomorphic
respectively metric incarnation of the same object. We will use the two
names accordingly.

\bigskip
{\bf \refstepcounter{theorem}\label{Notat} \thetheorem} ---
 Let $X$ be an irreducible symplectic manifold.
For later use we introduce the following notations:

The {\it K\"ahler cone}
$\kk_X\subset H^{1,1}(X)_\IR:=H^{1,1}(X)\cap H^2(X,\IR)$
is the open convex cone of all K\"ahler classes on $X$. 

If $\alpha\in H^{1,1}(X)_\IR$, then one defines
$$P(X):=(\sigma\cdot\IC\oplus\bar\sigma\cdot\IC)\cap H^2(X,\IR)
\phantom{XXX}{\rm and}\phantom{XXX}F(\alpha):=P(X)\oplus\alpha\cdot\IR,$$
where $\sigma$ is a holomorphic two-form spanning $H^0(X,\Omega_X^2)$. If
$\alpha\in\kk_X$ and the induced hyperk\"ahler structure is $(M,g,I,J,K)$,
where $X=(M,I)$, then $\sigma=\omega_J+i\omega_K$ (up to scalar
factors) and, therefore,
$P(X)=[\omega_J]\cdot\IR\oplus[\omega_K]\cdot\IR$ and
$F(\alpha)=[\omega_I]\cdot\IR\oplus[\omega_J]\cdot\IR\oplus[\omega_K]\cdot\IR$.
Hence, $F(\alpha)\subset H^2(X,\IR)=H^2(M,\IR)$ is independent of the complex
structure and depends only on the metric $g$ on $M$. Sometimes, $F(\alpha)$ is
called the {\it HK 3-space} associated with the hyperk\"ahler metric $g$.

\bigskip
{\bf \refstepcounter{theorem}\label{Cohomology} \thetheorem} ---
Let $X$ be an irreducible symplectic manifold. The following identities
are immediate consequences of the definition:
$$\begin{array}{rcl}
H^1(X,\ko_X)&=&0,~~ H^2(X,\ko_X)\cong\IC\\
H^0(X,\kt_X)&\cong&H^0(X,\Omega_X)=0\\
H^1(X,\kt_X)&\cong& H^1(X,\Omega_X).\\
\end{array}$$
In the following (\ref{hypcoh},\ref{irredcoh},\ref{pairing})
we state some results concerning the cohomological
structure of irreducible symplectic manifolds. Most of them are due to
Fujiki. For details and other results in this direction we refer to Enoki's
survey article \cite{Enoki}, to the original paper of Fujiki
\cite{Fujiki2} and to the more recent articles by Looijenga and Lunts \cite{LL}
and Verbitsky \cite{Verbitsky}.

\bigskip
{\bf \refstepcounter{theorem}\label{hypcoh} \thetheorem} ---
Let $(M,g)$ be a compact irreducible hyperk\"ahler manifold of real
dimension $4n$. Let $F$ denote
the associated HK $3$-space spanned by the three K\"ahler
forms $[\omega_I],[\omega_J],[\omega_K]$. The Lefschetz operator $L_\lambda:
=L_{\omega_\lambda}$ on the K\"ahler manifold $(M,\lambda,\omega_\lambda)$ for
$\lambda\in S^2$ acts on the cohomology $H^*(M,\IR)$ and
allows one to define the primitive cohomology $H^k((M,\lambda),\IR)_{\rm pr}:=
\ker(L_\lambda^{2n-k+1}:H^k(M,\IR)\to H^{4n-k+2}(M,\IR))$ ($k\leq2n$). Recall,
the
{\it Lefschetz decomposition} on the compact K\"ahler manifold
$(M,\lambda,\omega_\lambda)$ is the direct sum decomposition
$$H^k(M,\IR)=H^k((M,\lambda),\IR)=\bigoplus_{(k-\ell)\geq(k-2n)^+}L_\lambda^{k-\ell}(H^{2\ell-k}((M,\lambda),\IR)_{\rm pr})$$
and the {\it Hard Lefschetz Theorem} asserts that for $k\leq 2n$ 
$$L_\lambda^{2n-k}:H^k(M,\IR)\cong H^{4n-k}(M,\IR).$$
The latter implies
$$L_\lambda^{k-\ell}:H^{2\ell-k}((M,\lambda),\IR)_{\rm pr}\cong L_\lambda^{k-\ell}(H^{2\ell-k}((M,\lambda),\IR)_{\rm pr}).$$
Combining the statements for all $\lambda\in S^2$ one obtains the following
results \cite{Fujiki2}: Let $N^*\subset H^*(M,\IR)$ denote the subalgebra
generated by $F$ and let $H^*(M,\IR)_F:=\bigoplus H^k(M,\IR)_F=\bigoplus
\{\alpha\in H^k((M,\lambda),\IR)_{\rm pr}| ~{\rm for~all~}\lambda\in S^2\}$.
Then
$$H^k(M,\IR)=\bigoplus_{2(k-l)\geq(k-2n)^+} N^{k-\ell}H^\ell(M,\IR)_F$$
and for $k\leq n$
$$N^{k-\ell}\otimes_\IR H^\ell(M,\IR)_F\cong N^{k-\ell}H^\ell(M,\IR)_F.$$
The first statement in particular says
$$H^2(M,\IR)=F\oplus H^2(M,\IR)_F.$$
Note that for any $\lambda\in S^2$ the space $H^2(M,\IR)_F$ is contained
in $H^{1,1}((M,\lambda),\IR)$.
   
\bigskip
{\bf \refstepcounter{theorem}\label{irredcoh} \thetheorem} ---
Let $X$ be an irreducible symplectic manifold and let
$0\ne\sigma\in H^0(X,\Omega_X^2)$ be fixed. By the holonomy principle one
easily obtains (cf.\ \cite{Beauville1}):
$$H^0(X,\Omega_X^p)\cong\left\{\begin{array}{lcl}
0&&p\equiv1(2)\\
\Lambda^{p/2}\sigma\cdot\IC&&p\equiv0(2).\\
\end{array}\right.$$
Fujiki also proved holomorphic versions of the Lefschetz decomposition
and the Hard Lefschetz Theorem. They allow one
to compute other cohomology groups
of the form $H^q(X,\Omega_X^p)$ as follows: Let
$L_\sigma:H^q(X,\Omega^p_X)\to H^q(X,\Omega_X^{p+2})$ and
$L_{\bar\sigma}:H^q(X,\Omega_X^p)\to H^{q+2}(X,\Omega_X^p)$ be the
map given by the
cup-product with $\sigma$ and $\bar\sigma$, respectively, and let
$H^q(X,\Omega^p_X)_\sigma:= \ker(L_\sigma^{n-p+1})$ and
$H^q(X,\Omega_X^p)_{\bar\sigma}:=\ker(L_{\bar\sigma}^{n-q+1})$. Then
$$H^q(X,\Omega_X^p)
=\bigoplus_{(p-\ell)\geq(p-n)^+} L_\sigma^{p-\ell}H^q(X,\Omega_X^{2\ell-p})_\sigma
= \bigoplus_{(q-\ell)\geq(q-n)^+}L_{\bar\sigma}^{q-\ell}H^{2\ell-q}(X,\Omega_X^p)_{\bar\sigma}$$
and
$$\begin{array}{crcll}
L_\sigma^{n-p}:&H^q(X,\Omega_X^p)&\cong&H^q(X,\Omega_X^{2n-p})~~~~~&{\rm for}~p\leq n\\
L_{\bar\sigma}^{n-q}:&H^q(X,\Omega_X^p)&\cong&H^{2n-q}(X,\Omega_X^p)~~~~~&{\rm for}~q\leq n.\\
\end{array}$$
The following special cases will be used frequently
$$\begin{array}{crcl}
L_\sigma^{n-1}:&H^q(X,\Omega_X)&\cong&H^q(X,\Omega_X^{2n-1})~~~~~~{\rm for~all}~q\\
L_{\bar\sigma}^{n-1}:&H^1(X,\Omega_X^p)&\cong&H^{2n-1}(X,\Omega_X^p)~~~~~~{\rm for~all}~p.\\
\end{array}$$ 
The first isomorphism can be deduced using the isomorphism
$\sigma:\kt_X\to\Omega_X$ and
the perfect pairing $\Omega_X\otimes\Omega_X^{2n-1}\to K_X\cong\ko_X$:
The induced map $\Omega_X\cong\kt_X\cong(\Omega_X)\dual\cong\Omega_X^{2n-1}$
is just $L_\sigma^{n-1}$. The complex conjugate of the first isomorphism
gives the second. Combining both we get
$$L^{n-1}_{\sigma\bar\sigma}:=L_\sigma^{n-1}\circ L_{\bar\sigma}^{n-1}:H^1(X,\Omega_X)\cong H^{2n-1}(X,\Omega_X^{2n-1}).$$ 
In Sect.\ \ref{projectivity} we will also need a version of this isomorphism
on the level of smooth forms.

\bigskip
{\bf \refstepcounter{theorem}\label{pairing} \thetheorem} ---
The following natural pairings are perfect:
$$\begin{array}{ccccl}
H^1(X,\kt_X)&\otimes &H^{2n-1}(X,\Omega_X^{2n-1})&\to&H^{2n}(X,\Omega_X^{2n-2})\cong({\bar\sigma}^n\sigma^{n-1})\cdot\IC\\
H^1(X,\kt_X)&\otimes& H^1(X,\Omega_X)&\to&H^2(X,\ko_X)\cong\bar\sigma\cdot\IC\\
\end{array}$$

Using Serre duality it suffices to show that the induced homomorphisms
$H^1(X,\kt_X)\to\Hom(H^{2n-1}(X,\Omega_X^{2n-1}), H^{2n}(X,\Omega_X^{2n-2}))\cong\Hom(H^0(X,\Omega_X^2),H^1(X,\Omega_X))$ and
$H^1(X,\kt_X)\to \Hom(H^1(X,\Omega_X),H^2(X,\ko_X))\cong \Hom(H^{2n-2}(X,\Omega^{2n}_X),H^{2n-1}(X,\Omega_X^{2n-1}))$ are bijective.
Using the general form of the holomorphic
Hard Lefschetz Theorem above one finds that these maps are also given by
the cup-product. Then the perfectness of the first pairing follows
from the fact that $\sigma:\kt_X\to\Omega_X$ is bijective
and the proof of the second uses in addition that
the composition of $L^n_\sigma:H^1(X,\kt_X)\to H^1(X,\Omega_X^{2n-1})$
followed by
$L^{n-1}_{\bar\sigma}:H^1(X,\Omega^{2n-1}_X)\to H^{2n-1}(X,\Omega_X^{2n-1})$ is bijective.

Note that in particular for any $0\ne\alpha\in H^1(X,\Omega_X)$ (or
$\in H^{2n-1}(X,\Omega_X^{2n-1})$) the natural map
$\alpha:H^1(X,\kt_X)\to H^2(X,\ko_X)$ (respectively
$H^1(X,\kt_X)\to H^{2n}(X,\Omega_X^{2n-2})$)
is surjective.

We will also use the following fact:
Let $\beta\in H^{2n-1}(X,\Omega_X^{2n-1})$, then the cup-product with $\beta$
sends $v\in H^1(X,\kt_X)$ to $\beta\cdot v\in H^{2n}(X,\Omega_X^{2n-2})$.
Hence $(\beta\cdot v)\sigma$ can be integrated over $X$. On the other hand,
$\beta$ can be regarded as a linear form on $H^1(X,\Omega_X)$, and thus
can be evaluated on the image $\alpha$
of $v$ under the isomorphism $H^1(X,\kt_X)\cong H^1(X,\Omega_X)$.
Clearly, the two expressions
$\int_X(\beta\cdot v)\sigma$ and $\beta(\alpha)=\int_\beta\alpha$ agree.
A similar statement holds for $\beta\in H^1(X,\Omega_X)$. Here
we have $\int (\beta\cdot v)\sigma^n{\bar\sigma}^{n-1}=\int(\beta\alpha)(\sigma\bar\sigma)^{n-1}=c q_X(\beta,\alpha)$,
where $c$ is a positive number (for the definition of $q_X$ see
\ref{quadraticform}).

\bigskip
{\bf \refstepcounter{theorem}\label{quadraticform} \thetheorem} ---
Due to work of Beauville \cite{Beauville1} there exists a natural
quadratic form on the second cohomology of an irreducible symplectic
manifold generalizing the intersection pairing on a K3 surface. There are
several approaches towards this quadratic form
\cite{Beauville1,Fujiki2,Enoki,Bo3,Verbitsky} most
of them are intimately interwoven with the deformation theory
of such manifolds. Let us state the main facts.

Let X be an irreducible symplectic manifold and let
$\sigma\in H^0(X,\Omega^2_X)$ such that $\int(\sigma\bar\sigma)^n=1$.
Define a quadratic form on $H^2(X,\IR)$ by
$$f(\alpha)=\frac{n}{2}\int(\sigma\bar\sigma)^{n-1}\alpha^2+(1-n)\left(\int\sigma^{n-1}{\bar\sigma}^n\alpha\right)\cdot\left(\int\sigma^n{\bar\sigma}^{n-1}\alpha\right).$$
Writing $\alpha$ according to the Hodge decomposition as $\alpha=\lambda\sigma
+\beta+\bar\lambda\bar\sigma$, where $\beta$ is a $(1,1)$-form, then
$$f(\alpha)=\frac{n}{2}\int(\sigma\bar\sigma)^{n-1}\beta^2+\lambda\bar\lambda.$$
Clearly, with respect to this quadratic form $H^{1,1}(X)$ is orthogonal
to $H^{2,0}(X)\oplus H^{0,2}(X)$. Moreover, $f(\sigma)=0$ and
$f(\sigma+\bar\sigma)=1>0$. If $\alpha$ is a K\"ahler class on $X$ and
$F=F(\alpha)$,
then the decomposition $H^2(X,\IR)=F\oplus H^2(X,\IR)_F$ is orthogonal
with respect to $f$. The restriction of $f$ to $F$ only depends
on the underlying manifold $M$ and the hyperk\"ahler metric but not on
the complex structure.
One also has the following extremely useful formula
\beeq{extrform}
v(\alpha)^2f(\beta)=f(\alpha)\left((2n-1)v(\alpha)\int\alpha^{2n-2}\beta^2-(2n-2)\left(\int\alpha^{2n-1}\beta\right)^2\right)
\eneq
for any two classes $\alpha$ and $\beta$, where $v(\alpha)=\int\alpha^{2n}$.
Applying this formula to a K\"ahler class $\alpha$ and $\beta=\sigma+\bar\sigma$
yields $v(\alpha)=2f(\alpha)(2n-1)\int\alpha^{2n-2}(\sigma\bar\sigma)$.
By the Hodge-Riemann bilinear relations the integral on the right hand side
is positive. Hence: For any K\"ahler class $\alpha$ the quadratic form
$f$ restricted to $F(\alpha)$ is positive definite.
For $0\ne\beta\in H^2(X,\IR)_F$ the above formula shows
$v(\alpha)f(\beta)=(2n-1)f(\alpha)\int\alpha^{2n-2}\beta^2$. Thus
$f$ restricted to $H^2(X,\IR)_F$ is a positive
multiple of the standard Hodge-Riemann bilinear form and, therefore,
negative definite.
Yet another consequence of formula (\ref{extrform}) is the following: Since
$v(\sigma+\bar\sigma)>0$ and $f(\sigma+\bar\sigma)=1$, for any
$\beta\in H^2(X,\IQ)$ close to $\sigma+\bar\sigma$ one has
$f(\beta)>0$ and $v(\beta)\in\IQ$. Hence
$$f(\alpha)/f(\beta)=\frac{1}{v(\beta)^2}\left((2n-1)
v(\beta)\int\beta^{2n-2}\alpha^2-(2n-2)\left(\int\beta^{2n-1}\alpha\right)^2\right)\in\IQ$$
for all $\alpha\in H^2(X,\IQ)$. 

The upshot is: There exists a positive constant $c\in\IR$ such that
$q_X:=c\cdot f$ is a {\it primitive integral quadratic form} on
$H^2(X,\IZ)$ of index $(3,b_2(X)-3)$. For any K\"ahler class $\alpha$
the decomposition $F(\alpha)\oplus H^2(X,\IR)_{F(\alpha)}$ is orthogonal
with respect to $q_X$. Also note $q_X(\sigma)=0$ and
$q_X(\sigma+\bar\sigma)>0$. 

By means of the integral quadratic form $q_X$ one can establish
a close link between rational classes of dimension one and
those of codimension
one. Namely, if $c$ is the positive constant such that
$c\cdot f=q_X$, then $cL_{\sigma\bar\sigma}^{n-1}$
defines an isomorphism $H^{1,1}(X)_\IQ\cong H^{2n-1,2n-1}(X)_\IQ$.
(Note that this is quite similar to the Hard Lefschetz Theorem with
respect to a Hodge class.) Indeed, if $\alpha\in H^{1,1}(X)$ and
$\beta:=c L_{\sigma\bar\sigma}^{n-1}(\alpha)$, then
for any $\gamma\in H^2(X,\IQ)$ one has
$$\int\beta\gamma=\int\beta\gamma^{1,1}=c\int L_{\sigma\bar\sigma}^{n-1}(\alpha)\gamma^{1,1}=\frac{2}{n}q_X(\alpha,\gamma^{1,1})=\frac{2}{n}q_X(\alpha,\gamma).$$
Now, if $\alpha$ is rational then $q_X(\alpha,\gamma)\in\IQ$ and hence
$\int\beta\gamma\in\IQ$ for all $\gamma\in H^2(X,\IQ)$. Thus $\beta\in H^{4n-2}(X,\IQ)$. Conversely, if $\beta$ is rational, then $\int\beta\gamma\in\IQ$ and
therefore $q_X(\alpha,\gamma)\in\IQ$. Since $q_X$ is non-degenerate
and integral, this
shows $\alpha\in H^2(X,\IQ)$.

\bigskip
{\bf \refstepcounter{theorem}\label{posconedef} \thetheorem} ---
The {\it positive cone} $\kc_X$ is by definition the component
of $\{\alpha\in H^{1,1}(X)_\IR|q_X(\alpha)>0\}$ that contains the K\"ahler
cone $\kk_X$ (cf.\ \ref{Notat}). Note that one has a Hodge Index Theorem
with respect to $q_X$: A $(1,1)$-class $\beta$ which is orthogonal
to a K\"ahler class (with respect to $q_X$) satisfies $q_X(\beta)<0$
or is zero. In particular, if $\alpha\in \kk_X$ then
$q_X(\alpha,\,.\,)$ is positive on $\kc_X$. Moreover, if $\alpha\in\kk_X$ then
$\int\alpha^{2n-1}\beta>0$ for $\beta\in\kc_X$. Indeed,
if not then one could find $\beta\in\kc_X$ with $\int\alpha^{2n-1}\beta=0$,
i.e.\ $\beta\in H^2(X,\IR)_{F(\alpha)}$, which would imply
$q_X(\beta)<0$. This is absurd. 

\bigskip
{\bf \refstepcounter{theorem}\label{Todd} \thetheorem} ---
For the following we refer to Fujiki's paper \cite{Fujiki2} but we also
wish to draw the reader's attention to the more recent preprint
of Looijenga and Lunts \cite{LL}.

For any integral class $\alpha\in H^{2j}(X,\IZ)$ one has the form
of degree $2n-j$ that sends $\beta\in H^2(X,\IZ)$ to
$\int\alpha\beta^{2n-j}\in\IZ$. E.g.\ for $j=0$ and $\alpha=1$ this
is $v(\beta)=\int\beta^{2n}$. Fujiki shows that for any
$\alpha\in H^{4j}(X,\IZ)$ contained in the subalgebra generated
by the Chern classes of $X$ there exists a constant $c\in \IQ$ such that
\beeq{extrform2}
\int\alpha\beta^{2(n-j)}=cq_X(\beta)^{n-j}~~{\rm for~any}~
\beta\in H^2(X,\IQ).
\eneq
In particular, for $j=0$ and
$\alpha=1$ this yields $v(\beta)=\int\beta^{2n}=cq_X(\beta)^n$ and in this
case one has $c>0$. In fact, the result can be slightly generalized to all
classes $\alpha$ which are of type $(2j,2j)$ on all small deformations
of $X$. 

As an application of (\ref{extrform2})
one has that the Hirzebruch-Riemann-Roch formula on an
irreducible symplectic manifold takes the following form:
If $L$ is a line bundle on $X$ then
$$\chi(L)=\sum \frac{a_i}{(2i)!}q_X(c_1(L))^i,$$
where the $a_i$'s are constants only
depending on $X$. Indeed, $\chi(L)=\sum\int ch_{i}(L)td_{2n-i}(X)$ and
$(2i)! ch_{2i}(L)=c_1(L)^{2i}$. Since $td_i(X)$ is a
polynomial in the Chern classes, one has
$td_i(X)=0$ for $i\equiv 1(2)$ (use $c_i(X)=0$ for $i$ odd)
and $c_1(L)^{2i}td_{2n-2i}=a_iq_X(c_1(L))^i$.
This will be crucial in the proof of Theorem \ref{birat}.

Another application of the relation between $q_X(\alpha)$ and 
$v(\alpha)$ is the fact that if $X$ is irreducible symplectic,
$\alpha$ a K\"ahler class, and $Y\subset X$ an effective divisor,
then $q_X(\alpha,[Y])>0$. One possible proof goes as follows:
If $[Y]\in\kc_X$ then certainly $q_X(\alpha,[Y])>0$ by the Hodge Index Theorem.
Since $\int_Y\alpha^{2n-1}>0$ and hence $\int\alpha^{2n-1}(-[Y])<0$,
the case $-[Y]\in\kc_X$ can be excluded (cf.\ \ref{posconedef}).
Thus it remains
to verify the claim for $q_X([Y])\leq0$.
Of course, $q_X(\alpha,[Y])>0$ if and only if $q_X(\alpha+[Y])>
q_X(\alpha)+q_X([Y])$. Thus it suffices to show $q_X(\alpha+[Y])>
q_X(\alpha)$.
Replacing $\alpha$ by $k\beta$ with $k\gg0$ and $\beta$ fix, this is
equivalent to $\int (k\beta+[Y])^{2n}=c q_X(k\beta+[Y])^{n}>
c q_X(k\beta)^n=\int(k\beta)^{2n}$. Since $\int_Y\beta^{2n-1}>0$, this
follows immediately.

\bigskip
{\bf \refstepcounter{theorem}\label{deformationallg} \thetheorem} ---
A {\it deformation} of a compact manifold $X$ is a smooth proper holomorphic
map $\kx\to S$, where $S$ is an analytic space and the fibre over
a distinguished point $0\in S$ is isomorphic to $X$. We will say that
a certain property holds for the {\it generic} fibre, if for an open
(in the analytic topology) dense set $U\subset S$
and all $t\in U$ the fibre $\kx_t$ has this property. The property
holds for the {\it general} fibre if such a set $U$ exists that is the
complement of the union of
countably many nowhere dense closed (in the analytic topology)
subsets.

One knows that
for any compact K\"ahler manifold $X$ there exists a {\it semi-universal}
deformation $\kx\to \Def(X)$, where $\Def(X)$ is a germ of an analytic space and
the fibre $\kx_0$ over $0\in \Def(X)$ is isomorphic to $X$. The {\it Zariski
tangent space} of $\Def(X)$ is naturally isomorphic to $H^1(X,\kt_X)$.
If $H^0(X,\kt_X)=0$, i.e.\ if $X$ does not allow infinitesimal
automorphisms, then $\kx\to \Def(X)$ is universal, i.e.\ for any
deformation $\kx_S\to S$ of $X$ there exists a uniquely determined holomorphic
map $S\to \Def(X)$ such that $\kx_S\cong\kx\times_{\Def(X)}S$. By a result of
Tian \cite{Tian} and Todorov \cite{Todorov} the base space $\Def(X)$ is
smooth if $K_X\cong \ko_X$ (for an algebraic proof of this result we refer to
\cite{Kawamata} and \cite{R}, see also \cite{Bo2} for irreducible
symplectic manifolds). In this case one says that $X$ deforms
{\it unobstructed}.

Hence, if $X$ is an irreducible symplectic manifold then there exists a
universal deformation $\kx\to \Def(X)$ of $X$ with $\Def(X)$ smooth of dimension
$h^1(X,\kt_X)=h^1(X,\Omega_X)=h^{1,1}(X)$. Note that any small
deformation of an irreducible symplectic manifold is irreducible
symplectic (\cite{Beauville1}, see also \ref{exbydef}). Also note that the
universal deformation $\kx\to\Def(X)$ is in fact universal for any
fibre $\kx_t$ with $t$ close to $0$.

\bigskip
{\bf \refstepcounter{theorem}\label{deformationhk} \thetheorem} ---
Let $\alpha$ be a K\"ahler class on an irreducible symplectic
manifold $X$ and let $(M,g)$ be the underlying hyperk\"ahler manifold.
As briefly explained above, there exists a sphere $S^2\cong\IP^1$ of
complex structures on $M$ induced by the hyperk\"ahler metric $g$. 
This gives rise to a compact manifold $\kx(\alpha)$ and a smooth
holomorphic map 
$$\kx(\alpha)\to \IP^1$$
such that for any $\lambda\in S^2\cong\IP^1$
the fibre over $\lambda$ is isomorphic to $(M,\lambda)$.
In particular, $X$ occurs as the fibre over $I$. The space
$\kx(\alpha)$ together with the projection is called the
{\it twistor space} of the hyperk\"ahler manifold $(M,\lambda)$.
For more details and the precise relation between the metric $g$ and the
twistor space see \cite{HKR}.

The twistor space $\kx(\alpha)\to \IP^1$ induces a non-trivial (local) map
$\IP^1\to \Def(X)$ and hence a one-dimensional subspace of
$H^1(X,\kt_X)\cong H^1(X,\Omega_X)$. Fujiki \cite{Fujiki2} shows
that this subspace is spanned by $\alpha$. In other words, the Kodaira-Spencer
class of the twistor space of $X$ with respect to the K\"ahler class $\alpha$ 
is identified with $\alpha$ under the isomorphism $H^1(X,\kt_X)\cong
H^1(X,\Omega_X)$ induced by the holomorphic two-form $\sigma$ on $X$.
Note that one can also prove the unobstructedness of $X$ along this line:
Any class in $\kk_X\subset H^1(X,\Omega_X)$ can be realized as a
Kodaira-Spencer class of a smooth curve and $\kk_X$ is open in the real part
$H^{1,1}(X)_\IR$ of $H^{1,1}(X)$ (cf.\ \cite{Fujiki2}).

\bigskip
{\bf \refstepcounter{theorem}\label{deformationlb} \thetheorem} ---
At several points in this article we will have to deal with joint
deformations of an irreducible symplectic manifold
$X$ together with a line bundle $L$ on $X$. Again, there
exists a universal deformation $(\kx,\kl)\to \Def(X,L)$, i.e.\
a deformation $\kx\to \Def(X,L)$ of $\kx_0=X$ and a line bundle $\kl$
on $\kx$ with $\kl_0:=\kl|_{\kx_0}\cong L$ such that any other deformation
$(\kx_S,\kl_S)\to S$ of this sort is isomorphic to the pull-back of
$(\kx,\kl)$ via a uniquely determined map $S\to \Def(X,L)$.
The Zariski tangent space of
$\Def(X,L)$ is canonically isomorphic to $H^1(X,\kd(L))$, where $\kd(L)$
is the sheaf of differential operators on $L$ of order $\leq1$. Using
the symbol sequence and $H^1(X,\ko_X)=0$ one shows
$H^1(X,\kd(L))=\ker(c_1(L):H^1(X,\kt_X)\to H^2(X,\ko_X))$.
By \ref{pairing} the contraction
$c_1(L):H^1(X,\kt_X)\to H^2(X,\ko_X)=\bar\sigma\cdot\IC$ is surjective
if $0\ne c_1(L)\in H^1(X,\Omega_X)$. Moreover, in this case the space
$\Def(X,L)$ is a smooth hypersurface of $\Def(X)$. Another consequence
of \ref{pairing} is that for two line bundles $L$ and $M$ with $c_1(L)$
and $c_1(M)$ linearly independent the two hypersurfaces
$\Def(X,L)$ and $\Def(X,M)$ intersect transversally.

\bigskip
{\bf \refstepcounter{theorem}\label{periodmapdef} \thetheorem} ---
Let $\Gamma$ be a lattice of index $(3,b-3)$. By $q_\Gamma$ we denote its
quadratic form. A {\it marked} irreducible symplectic manifold is
a tuple $(X,\varphi)$ consisting of an irreducible symplectic manifold
$X$ and an isomorphism $\varphi:H^2(X,\IZ)\cong \Gamma$ compatible
with $q_X$ and $q_\Gamma$. The {\it period} of $(X,\varphi)$ is by definition
the one-dimensional subspace $\varphi(H^{2,0}(X))\subset \Gamma_\IC$
considered as a point in the projective space $\IP(\Gamma_\IC)$.
If $\kx\to \Def(X)$ is the universal deformation of $\kx_0=X$, then
a marking $\varphi$ of $X$ naturally defines markings $\varphi_t$ of all the
fibres $\kx_t$. Thus we can define the {\it period map}
$$\kp:\Def(X)\to \IP(\Gamma_\IC)$$
as the map that takes $t$ to the period of $(\kx_t,\varphi_t)$. Note
$\kp$ is holomorphic. Its tangent map is given by the contraction
$$H^1(X,\kt_X)\to \Hom(H^{2,0}(X), H^{1,1}(X))\subset \Hom(H^{2,0}(X), H^{2}(X,\IC)/H^{2,0}(X)).$$
By \ref{quadraticform} the holomorphic two-form $\sigma$ on $X$ satisfies
$q_X(\sigma)=0$ and $q_X(\sigma+\bar\sigma)>0$.
Hence the image of $\kp$ is contained in the {\it period domain}
$Q\subset \IP(\Gamma_\IC)$ defined as $\{x\in\IP(\Gamma_\IC)|q_\Gamma(x)=0,~~
q_\Gamma(x+\bar x)>0\}$, which is an open (in the analytic topology)
subset of the non-singular quadric defined by $q_\Gamma$.

Beauville proved in \cite{Beauville1} the Local Torelli Theorem:
{\it For any marked irreducible symplectic manifold $(X,\varphi)$ the
period map $\kp:\Def(X)\to Q$ is a local isomorphism.}

\bigskip
{\bf \refstepcounter{theorem}\label{hypersurfaces} \thetheorem} ---
Let $(X,\varphi)$ be a marked irreducible symplectic manifold. For any $\alpha
\in H^2(X,\IR)$ we define 
$$S_\alpha:=\{t\in \Def(X)|q_\Gamma(\varphi(\alpha),\kp(t))=0\},$$
i.e.\ $S_\alpha$ is the pull-back of the hyperplane defined by
$q_\Gamma(\varphi(\alpha),\,.\,)$. By the properties of $q_X$, the set
$S_\alpha$ is easily identified as the set of points $t\in \Def(X)$ such that
$\alpha$ is of type $(1,1)$ on $\kx_t$. Analogously, if
$\alpha\in H^{4n-2}(X,\IR)$, which can be considered as a linear form
on $H^2(X,\IR)$, then
$$S_\alpha:=\{t\in \Def(X)|\alpha(\varphi^{-1}(\kp(t)))=0\}.$$
Using the holomorphic Lefschetz theorem \ref{irredcoh}
one finds that $S_\alpha$ in this situation
is the set of $t\in \Def(X)$ such that $\alpha$ is of type
$(2n-1,2n-1)$ on the fibre $\kx_t$. 

Using the perfectness of the natural pairings \ref{pairing} one
proves the following results \cite{Beauville1,Fujiki2}:
\begin{itemize}
\item[--] For $0\ne\alpha\in H^2(X,\IR)$ or
$0\ne\alpha\in H^{4n-2}(X,\IR)$ the set
$S_\alpha$ is a smooth possibly empty
hypersurface of $\Def(X)$. 
\item[--] If $0\in S_\alpha$, i.e.\
$\alpha\in H^{1,1}(X)_\IR$ respectively $\alpha\in H^{2n-1,2n-1}(X)_\IR$,
the tangent space $T_0S_\alpha$ of $S_\alpha$
is the kernel of the surjection $\alpha:H^1(X,\kt_X)\to H^2(X,\ko_X)$
respectively $\alpha:H^1(X,\kt_X)\to H^{2n}(X,\Omega_X^{2n-2})$.
\item[--] If $\alpha,\alpha'\in H^{1,1}(X)$ are linearly independent, then
$S_\alpha$ and $S_{\alpha'}$ intersect transversally in $0\in \Def(X)$.
The analogous statement holds true
for linearly independent classes $\alpha,\alpha'\in H^{2n-1,2n-1}(X)$.
\item[--] If $\alpha\in H^2(X,\IZ)$, then there exists a line bundle $\kl$ on
$\kx|_{S_\alpha}$ such that $c_1(\kl_t)=\alpha$ for all $t\in S_\alpha$.
Moreover, if $\alpha=c_1(L)$, then $S_\alpha=\Def(X,L)$.
\end{itemize}

\bigskip
{\bf \refstepcounter{theorem}\label{twistor2} \thetheorem} ---
In \ref{deformationhk} we discussed the twistor
space $\kx(\alpha)\to \IP^1$ of an irreducible symplectic
manifold $X$ endowed with a K\"ahler class $\alpha$. The base of
the twistor space can be identified via the period map as follows:
Recall,
$F(\alpha)_\IC=\sigma\cdot\IC\oplus\bar\sigma\cdot\IC\oplus\alpha\cdot\IC$.
If $\varphi$ is a marking of $X$, then
the period map $\kp:\IP^1\to \IP(\Gamma_\IC)$ defines an isomorphism
of $\IP^1$ with $Q\cap\IP(\varphi(F(\alpha)_\IC))$. Thus, $\IP^1$,
as the base space of the twistor space, is naturally identified
with the quadric $T(\alpha)\subset \IP(F(\alpha)_\IC)$ defined
by $q_X(\beta)=0$. In the sequel, we will denote the twistor
space by $\kx(\alpha)\to T(\alpha)$. This can be slightly generalized as
follows. If $\alpha$ is just any class of type $(1,1)$ on $X$ we
define $T(\alpha):=\Def(X)\cap\kp^{-1}(\IP(\varphi(F(\alpha)_\IC)))$
and $\kx(\alpha)$ as the restriction of the universal family to $T(\alpha)$.
Note that for a K\"ahler class $\alpha$ the space
$T(\alpha)$ means the complete base $\IP^1$ of the
twistor space, but in general it is only defined as a closed subset of
$\Def(X)$.

Fujiki \cite{Fujiki1}
proved the following very useful result: {\it If $\alpha$ is a K\"ahler
class and $\kx(\alpha)\to T(\alpha)$ is the associated twistor space, then
the general fibre  contains neither effective divisors nor curves.}

Indeed, the tangent space of the twistor space is spanned by the image $v$
of $\alpha$ under the isomorphism $H^1(X,\Omega_X)\cong H^1(X,\kt_X)$.
On the other hand, if $D\subset X$ is an effective divisor and $C\subset X$
is a curve, then the tangent space of the hypersurface $S_{[D]}$ (respectively
$S_{[C]}$) is the kernel of the map
$H^1(X,\kt_X)\to H^2(X,\ko_X)$ (respectively
$H^1(X,\kt_X)\to H^{2n}(X,\Omega_X^{2n-2})$).
Since up to scalar factors
$\int_X ([D]\cdot v) \sigma^{n}{\bar\sigma}^{n-2}=q_X(\alpha,[D])>0$
and $\int_X ([C]\cdot v)\sigma=\int_C\alpha\ne0$
(use \ref{Todd} resp.\ \ref{pairing} and that $\alpha$ is K\"ahler),
the hypersurfaces $S_{[D]}$ and $S_{[C]}$ meet $T(\alpha)$ locally only
in $0$. Or equivalently, neither $D$ nor $C$ deforms in the twistor space.
Since there are only countably many such classes $[D]$ and $[C]$, this proves
the claim.  

\bigskip
{\bf \refstepcounter{theorem}\label{moduli} \thetheorem} ---
Last but not least we fix the notation for the moduli space
of marked irreducible symplectic manifolds. Let $\Gamma$ be
a lattice of index $(3,b-3)$, where $b\geq3$. Then
$ \gM_\Gamma=\{(X,\varphi)\}/\sim$, where $(X,\varphi)$ is a marked
irreducible symplectic manifold (we usually also fix the dimension
$2n$) and $(X,\varphi)\sim(X',\varphi')$ if and only if there exists
an isomorphism $f:X\cong X'$ such that $f^*=\pm(\varphi^{-1}\circ\varphi')$.
The Local Torelli Theorem (\cite{Beauville1}, \ref{periodmapdef}) allows
one to patch the local charts $\Def(X)$. Thus $\gM_\Gamma$
carries the structure of a non-separated (i.e.\ non-Hausdorff)
complex manifold. The period map can
be considered as a holomorphic map
$\kp:\gM_\Gamma\to Q\subset\IP(\Gamma_\IC)$.


\section{Examples}\label{examples}

For the reader's convenience we collect the known
examples of irreducible symplectic manifolds. In all cases the verification
is reduced to the conditions of Definition \ref{IrredSympl}
and not to an explicit
construction of an irreducible hyperk\"ahler metric. Explicit
examples of irreducible hyperk\"ahler metrics on compact manifolds
would be highly desirable.

\bigskip
{\bf \refstepcounter{theorem}\label{K3} \thetheorem. K3 surfaces.}
{\it A complex manifold of dimension two
is irreducible symplectic if and only if it is a K3 surface.}

By definition a K3 surface is a compact connected surface with trivial
canonical bundle and vanishing first Betti number. One can show
that any K3 surface is deformation equivalent to a smooth
quartic hypersurface in $\IP^3$ and, therefore, simply connected.
That a K3 surface
is K\"ahler is due to Siu \cite{Siu} (see also \cite{Periodes}).
Examples of Guan \cite{Guan} show that in higher dimensions
not every simply connected holomorph-symplectic manifold is K\"ahler.

\bigskip
{\bf \refstepcounter{theorem}\label{HilbK3} \thetheorem.
Hilbert schemes of K3 surfaces.}
{\it If $S$ is a K3 surface, then
$\Hilb^n(S)$ is irreducible symplectic (cf.\ \cite{Beauville1})}

By the Hilbert scheme $\Hilb^n(S)$
we mean the Douady space of zero-dimensional subspaces $(Z,\ko_Z)$ of $S$
of length $\dim_\IC\ko_Z=n$. Strictly speaking, $\Hilb^n(S)$ is a scheme only
if $S$ is algebraic. In general, it is just a complex space.
Using that $S$ is smooth, compact, connected, and of dimension two,
one shows that $\Hilb^n(S)$ is a smooth compact connected manifold
of dimension $2n$. By results of Varouchas \cite{Varouchas} the Hilbert scheme
is K\"ahler if the underlying surface is K\"ahler which is the
case for K3 surfaces. Beauville then concluded that for
any K3 surface $S$ the Hilbert scheme $\Hilb^n(S)$ is irreducible symplectic
by showing that $\Hilb^n(S)$ admits a unique (up to scalars)
everywhere non-degenerate holomorphic two-form and that it is simply
connected. For $n=2$ this result was also obtained by Fujiki. It is interesting
to note that for $n>1$ one has $b_2(\Hilb^n(S))=23$. Moreover,
the second cohomology $H^2(\Hilb^n(S),\IZ)$ endowed with the natural
quadratic form $q_X$ (cf.\ \ref{quadraticform}) is isomorphic to the lattice
$H^2(S,\IZ)\oplus (-2(n-1)\cdot\IZ)$.

\bigskip
{\bf \refstepcounter{theorem}\label{Kummer} \thetheorem.
Generalized Kummer varieties.}
{\it If $A$ is a two-dimensional torus,
then
${\rm K}^{n+1}(A)$ is irreducible symplectic (cf.\ \cite{Beauville1}).}

The generalized
Kummer variety ${\rm K}^{n+1}(A)$ is by definition
the fibre over $0\in A$ of the
natural morphism $\Hilb^{n+1}(A)\to S^{n+1}(A)\rpfeil{5}{\Sigma}A$,
where $\Sigma$ is the summation and $0\in A$ is the zero-point of the torus. 
$\Hilb^{n+1}(A)$ itself also admits an everywhere non-degenerate two-form,
but neither is this two-form unique nor is $\Hilb^{n+1}(A)$ simply connected.
That both conditions are satisfied for ${\rm K}^{n+1}(A)$
was shown by Beauville. That
${\rm K}^{n+1}(A)$ is K\"ahler follows again from the results in \cite{Varouchas}.
The second Betti number of ${\rm K}^{n+1}(A)$ is $7$ (cf.\ \cite{Beauville1}).

I usually refer to the examples provided by the Hilbert schemes of K3
surfaces and by the generalized Kummer varieties as the two standard series
of examples of irreducible symplectic manifolds. Note that by means of these
examples we have in any real dimension $4n$ at least two different compact real
manifolds admitting irreducible hyperk\"ahler metrics. That they are not
diffeomorphic (in fact, not even homeomorphic) follows from a comparison
of their second Betti numbers.

\bigskip
{\bf \refstepcounter{theorem}\label{exbydef} \thetheorem.
New examples by deformation.} {\it Any small deformation of an irreducible
symplectic manifold is irreducible symplectic (cf.\ \cite{Beauville1}).}

The stability results in \cite{KSIII} say that any small deformation of
a compact K\"ahler manifold is again K\"ahler. Since the Hodge number $h^{2,0}$
is constant in families of compact K\"ahler manifolds, any
small deformation of an irreducible symplectic manifold
admits a unique non-trivial non-degenerate two-form
which is everywhere non-degenerate.
In fact, using the splitting theorem one can show that any K\"ahler deformation
of an irreducible symplectic manifold is again irreducible symplectic.
For the details see \cite{Beauville1}. It seems to be an open problem if any
deformation of an irreducible symplectic manifold is K\"ahler and,
therefore, irreducible symplectic.
Deforming the Hilbert scheme $\Hilb^n(S)$ of a K3 surface or the generalized
Kummer variety ${\rm K}^{n+1}(A)$ provides new examples of
irreducible symplectic manifolds. Indeed, for $n>1$ one has
$\dim \Def(\Hilb^n(S))=h^{1,1}(\Hilb^n(S))=h^{1,1}(S)+1=\dim \Def(S)+1(=21)$
and $\dim \Def({\rm K}^{n+1}(A))=h^{1,1}({\rm K}^{n+1}(A))=h^{1,1}(A)+1=\dim
\Def(A)+1(=5)$.
Thus one can think of the deformations of $\Hilb^n(S)$ that remain
Hilbert schemes as a hypersurface in the full deformation space of
$\Hilb^n(S)$. An analogous result holds true for Kummer varieties.

\bigskip
{\bf \refstepcounter{theorem}\label{exbybir} \thetheorem.
New examples by birational transformation.} {\it If a compact K\"ahler
manifold admits an everywhere non-degenerate two-form and is birational
to an irreducible symplectic manifold, then it is irreducible
symplectic as well.}

This follows easily from the observation that $\pi_1(X)=\pi_1(Y)$
and $h^{2,0}(X)=h^{2,0}(Y)$ for two birational compact manifolds $X$ and
$Y$. Can one drop the assumption that the manifold is K\"ahler?
Or equivalently, is any compact manifold that is birational to an
irreducible symplectic manifold and that admits an everywhere non-degenerate
two-form, automatically K\"ahler? Note also the following

\begin{lemma}\label{hodgeunderbir}---
If $X$ and $X'$ are birational irreducible symplectic
manifolds, then there exists a natural isomorphism $H^2(X,\IZ)\cong H^2(X',\IZ)$
compatible with the Hodge structures and the quadratic forms
$q_X$ and $q_{X'}$ \cite{Enoki, Mu2,OG}.
\end{lemma}

\prf We provide three slightly different descriptions of this
isomorphism. First, let us fix maximal open subsets $U\subset X$ and
$U'\subset X'$ such that $U\cong U'$ and $\codim(X\setminus U),\codim (X'\setminus U')\geq2$ (see the remarks after \ref{assump}).
Then one defines the isomorphism as the composition
$H^2(X,\IZ)\cong H^2(U,\IZ)\cong H^2(U',\IZ)\cong H^2(X',\IZ)$.
Since these isomorphisms commute with
$H^0(X,\Omega_X^2)\cong H^0(U,\Omega_U^2)\cong H^0(U',\Omega_{U'}^2)
\cong H^0(X',\Omega_{X'}^2)$, the isomorphism $H^2(X,\IZ)\cong H^2(X',\IZ)$
is compatible with the Hodge structures.

The second description is in terms of the closure of the graph
$Z\subset X\times X'$ of the birational map $X- - \to X'$.
Then $[Z]_*:H^2(X,\IZ)\to H^2(X',\IZ)$, defined by $\alpha\mapsto
p'_*([Z].p^*\alpha)$, equals the above isomorphism. We also write
$[Z]_*$ for the map $\beta\mapsto p_*([Z].{p'}^*\beta)\in H^2(X,\IZ)$.
Using the first description, one finds $[Z]_*\circ [Z]_*=\id$.
At one point in the discussion later on we will need
$q{_X'}([Z]_*\alpha,\beta)=q_X(\alpha,[Z]_*\beta)$, which follows
easily from the compatibility with the quadratic form, proved below,
and $[Z]_*\circ [Z]_*=\id$.

Yet another description goes as follows: Let $\tilde Z\to Z$
be a desingularization. Then
$H^2(\tilde Z,\IZ)\cong H^2(X,\IZ)\oplus\bigoplus_i[E_i]\cdot\IZ$ and
$H^2(\tilde Z,\IZ)\cong H^2(X',\IZ)\oplus\bigoplus_i[E_i]\cdot\IZ$ , where
the $E_i$'s are the exceptional divisors of $\tilde Z\to X$ (or,
equivalently, of $\tilde Z\to X'$).
The isomorphism $H^2(X,\IZ)\cong H^2(X',\IZ)$ is then given as
the composition of the natural
inclusion followed by the projection.

In either of the three descriptions one sees that the isomorphism
maps the class of an  effective divisor in $X$ to the class of an effective
divisor in $X'$.

To see the compatibility with the quadratic forms $q_X$ and $q_{X'}$, let 
$\sigma$ and $\sigma'$ be two-forms on $X$ and $X'$, respectively, with
$\int_X(\sigma\bar\sigma)^{n}=1=\int_{X'}(\sigma'{\bar\sigma}')^n$.
Their pull-backs to $\tilde Z$ coincide. Using the description
of $q_X$ and $q_{X'}$
given in \ref{quadraticform}, one finds that both can be defined on $\tilde Z$
and that via $[Z]_*$ they coincide if $\sigma^{n-1}|_{E_i}=0$.
The latter equality now follows from the observation
that the fibers (which are of positive dimension)
of the  morphisms $E_i\to X$ and $E_i\to X'$ are different and, hence,
the rank of $\sigma$ (or $\sigma'$) on $E_i$ drops at least by two.
\qed

\bigskip
There exist non-trivial birational transformations: Mukai introduced
the notion of {\it elementary transformations} \cite{Mu1}:
Let $X$ be an irreducible symplectic manifold and assume that 
a smooth $\IP^m$-bundle $P:=\IP(F)\rpfeil{5}{\phi}Y$ can be embedded
into $X$ as a submanifold of codimension $m$.
Then $\kn_{P/X}\cong\Omega_\phi$ and, hence, the exceptional divisor
of the blow-up $\tilde X\to X$ of $X$ along $P$ is isomorphic
to $\IP(\Omega_\phi)$. Regarding $\IP(\Omega_\phi)$ as the incidence
variety in $\IP(F)\times_Y\IP(F\dual)$ yields another projection
$\IP(\Omega_\phi)\to\IP(F\dual)$. Using the blow-down
criterion of Fujiki and Nakano one extends this projection to
a blow-down $\tilde X\to X'
(:=elm_P(X))$. It is not hard to see that $X'$ also admits a unique everywhere
non-degenerate two-form and that $X'$ is simply connected.
If $X'$ is K\"ahler (is this always true?)
then it is also irreducible symplectic.

\begin{example}\label{bevdeb}---
First, let us recall an example of Beauville \cite{Beauville2}: If $S\subset
\IP^3$ is a smooth quartic hypersurface, then the generic line
$\ell\subset \IP^3$ meets $S$ in four distinct points. Thus one defines
a rational map $\Hilb^2(S)- - \to \Hilb^2(S)$ by sending
$[Z]\in \Hilb^2(S)$ to $(\ell_Z\cap S)\setminus\{Z\}$, where $\ell_Z$ is the
uniquely defined line through $Z$. If $S$ does not contain any line then
this map extends to an automorphism of $\Hilb^2(S)$. On the other hand, 
if $S$ contains $k$ lines $\ell_1,\ldots,\ell_k$, then
$\Hilb^2(S)- - \to \Hilb^2(S)$ is the elementary transformation
along the naturally embedded planes $\IP^2\cong
\Hilb^2(\ell_i)\subset \Hilb^2(S)$. But, not any elementary transformation
of an irreducible symplectic manifold $X$ is again isomorphic to $X$!
One should rather think of Beauville's example as an exception. However,
it is usually not easy to show that $X'=elm_P(X)$ is not isomorphic
to $X$. One explicit example was given by Debarre in \cite{Debarre}:

Let $S$ be a K3 surface such that $\Pic(S)=\ko(C)\cdot\IZ$, where $C$
is a $(-2)$-curve. The Hilbert scheme $X:=\Hilb^n(S)$ contains the projective
space $\Hilb^n(C)\cong S^n(C)\cong\IP^n$ as above. Let $X'$ be the
elementary transformation of $X$ along this projective space. If
$X$ and $X'$ were isomorphic, then this would yield
a birational automorphism of $X$. Debarre then shows that,
due to the special shape of the surface $S$,
the induced map on the second cohomology would be trivial and
that this implies that the birational automorphism can be extended
to an isomorphism, which is absurd. Due to \ref{hodgeunderbir} the two
manifolds $X$ and $X'$ have isomorphic periods, but in order
to obtain an honest counterexample to the Global Torelli Theorem
one has in addition to ensure that $X'$ is K\"ahler.
This is not clear in general, but if $S$ is close to an algebraic
K3 surface $S_0$ as in the example above,
then $X'$ is close to $elm_{\IP^n}(\Hilb(S_0))\cong \Hilb^n(S_0)$.
But any small deformation of a K\"ahler manifold is K\"ahler; hence
there are examples for which $X'$ is K\"ahler. Note that in Debarre's example
the manifolds are not projective. It would be interesting to find
also a projective counterexample to the Global Torelli Theorem.
Good candidates are the moduli space of stable sheaves on
a K3 surface and the Hilbert scheme of the surface, which in some
cases are birational, but most likely not isomorphic (cf.\ \ref{modK3}).
\end{example}

The following examples are related to one of the standard series by
deformation or birational transformation as described above.
In particular, the proof that they are irreducible
symplectic is reduced to the proof of this fact for either
a Hilbert scheme or a Kummer variety.

\bigskip

{\bf \refstepcounter{theorem}\label{fano} \thetheorem.
Fano varieties of cubics \cite{BD}.} Let $Y\subset \IP^5$ be a smooth
cubic hypersurface and let $X:=F(Y)$ be the Fano variety of lines on $Y$.
Then $X$ is smooth and projective of dimension $4$. All cubics are deformation
equivalent and, hence, so are the corresponding Fano varieties.
Beauville and Donagi showed that for a special cubic $Y$ the Fano variety
$X=F(Y)$ is isomorphic to the Hilbert scheme $\Hilb^2(S)$ of a special
K3 surface of degree $14$ in $\IP^8$. Hence, for an arbitrary cubic $Y$ the
Fano variety $X=F(Y)$ is a deformation of $\Hilb^2(S)$ and, therefore,
irreducible symplectic. The map $Y\mapsto F(Y)$ identifies $\Def(Y)$
with a hypersurface of $\Def(X)$ parametrizing algebraic deformations
of $\Hilb^2(S)$. It would be interesting to understand the space of algebraic
deformations of $\Hilb^n(S)$ for other K3 surfaces.

\bigskip
{\bf \refstepcounter{theorem}\label{reljac} \thetheorem.
Relative Jacobians \cite{Markushevich}.} Markushevich constructed an
explicit example of a projective irreducible symplectic manifold
which is completely integrable: Let $\pi:S\to \IP^2$ be a generic double
cover ramified along a sextic. Then $S$ is a K3 surface. The dual space
${\IP^2}\dual$ can be regarded as the base space of the family of hyperelliptic
curves (of genus $2$) of the form $\pi^{-1}(\ell)$, where
$\ell\subset\IP^2$ is a line. Then the compactified relative Jacobian
$X\to{\IP^2}\dual$ of this family of curves is smooth projective and admits
an everywhere non-degenerate two-form.
Since the map that sends $[Z]\in \Hilb^2(S)$ to the divisor $(Z)$
on the curve $\pi^{-1}(\ell_{\pi(Z)})$,
where $\ell_{\pi(Z)}$ is the line through $\pi(Z)$, defines a birational
map $\Hilb^2(S)- -\to X$, the manifold $X$ is irreducible symplectic.
In fact, $X$ can be seen as an example of a moduli space
of simple sheaves (or rather stable sheaves with pure support
as treated below). Similar examples were considered by Beauville and Mukai.
But, as it is not surprising, in all cases the resulting spaces
are birational to some Hilbert scheme of the underlying K3 surfaces. 

\bigskip
{\bf \refstepcounter{theorem}\label{modK3} \thetheorem.
Moduli spaces of sheaves on K3 surfaces.} Let $S$ be a K3 surface.
Consider $M_H(v)$ -- the moduli space of sheaves with primitive Mukai
vector $v$ which are stable with respect to a generic polarization $H$
(see \cite{OG} or \cite{HL} for the notation).
This space is a projective variety and by the general smoothness
criterion it is non-singular. Note that also $\Hilb^n(S)$ can be
considered as such a moduli space, namely the moduli space of
stable rank one sheaves.
Mukai \cite{Mu1} constructed on $M_H(v)$ an everywhere
non-degenerate two-form. Later it was shown that $M_H(v)$ is indeed
irreducible symplectic (in \cite{GottscheHuybrechts} for the rank two case
and in \cite{OG} for arbitrary rank). The proof in both cases goes roughly
as follows: Any smooth deformation of $(S,H)$ induces a deformation
of $M_H(v)$ which is smooth as long as $H$ stays $v$-generic. Using this
one reduces the proof to the case of a very special K3 surface (e.g.\
an elliptic surface in \cite{OG}), where
one can show that the moduli space is birational to the Hilbert scheme
of the same dimension. Hence we have: {\it The moduli spaces $M_H(v)$ are
irreducible symplectic manifolds whenever $v$ is primitive and $H$ is
$v$-generic.}

\section{Projectivity}\label{projectivity}

Let $X$ be an irreducible symplectic manifold. The goal of this section is
to show that the second cohomology $H^2(X,\IZ)$ endowed with the natural
weight-two Hodge structure and the quadratic form $q_X$ (cf.\
\ref{quadraticform})
determines whether $X$ is a projective variety. The following list
collects known projectivity criteria that either motivate the main result
(Theorem \ref{proj}) or are essential for its proof.

\bigskip
\refstepcounter{theorem}\label{Moish}{\bf \thetheorem} ---
Let $X$ be a compact complex manifold.
Then $X$ is projective if and only if $X$ is K\"ahler and Moishezon
(cf.\ \cite{Moishezon}, \cite{Peternell}).

\refstepcounter{theorem}\label{Kod}{\bf \thetheorem} ---
Let $X$ be a compact complex manifold.
Then $X$ is projective if and only if $X$ admits a K\"ahler form
$\omega$ such that its cohomology class
$[\omega]\in H^2(X,\IC)$ is integral (cf.\ \cite{Kodaira2},
\cite{GriffithsHarris}).

\refstepcounter{theorem}\label{Chow}{\bf \thetheorem} ---
Let $X$ be a compact complex surface.
Then $X$ is projective if and only if $X$ is Moishezon (cf.\ \cite{Kodaira1},
\cite{BPV}).

\refstepcounter{theorem}\label{possurf}{\bf \thetheorem} ---
Let $X$ be a compact complex
surface. Then $X$ is projective if and only if there exists a line bundle $L$
on $X$ such that $\int c_1^2(L)>0$ (cf.\ \cite{BPV}).

\bigskip
It is easy to see that a surface $X$ that admits a line bundle $L$ with
$\int c_1^2(L)>0$ is Moishezon. Indeed, the Hirzebruch-Riemann-Roch formula
shows that either $h^0(X,L^m)\sim m^2$ or $h^0(X,L^{-m}\otimes K_X)\sim m^2$.
For manifolds of dimension $>2$ the existence of a big line bundle $L$, i.e.\
a line bundle with $\int c_1(L)^{\dim(X)}>0$, is not
sufficient to conclude that $X$ is Moishezon,
for the Hirzebruch-Riemann-Roch formula only gives $\sum h^{2i}(X,L^m)\sim
m^{\dim(X)}$.

Let us begin with a result due to Fujiki \cite{Fujiki1}
(see also \cite{Campana}), which in particular shows that projective
irreducible symplectic manifolds are dense in the moduli space.
Some of the techniques used in later chapters are based on the proof of
this result. Thus, we decided to reproduce it here.

\begin{theorem}\label{fujiki}{\rm \cite{Fujiki1}}---
Let $X$ be an irreducible symplectic manifold, let $0\in S\subset
\Def(X)$ be a smooth analytic subset of positive dimension
and let $\kx\to S$ denote the restriction of the Kuranishi family to $S$.
Then any open neighbourhood $U\subset S$ of $0\in S$ contains
a point $t\ne0$ such that the fibre $\kx_t$ over $t$ is projective.
\end{theorem}

\prf We may assume that $S$ is one-dimensional. Then, its tangent space at
zero $T_0S$ is a line in $T_0\Def(X)=H^1(X,\kt_X)$ spanned by, say,
$0\ne v\in H^1(X,\kt_X)$. By \ref{pairing} the induced map
$\tilde v:H^1(X,\Omega_X)\to H^2(X,\ko_X)\cong\IC$ is surjective. Since the
K\"ahler cone $\kk_X$ is an open subset of
$H^{1,1}(X)_\IR= H^{1,1}(X)\cap H^2(X,\IR)$, which spans $H^1(X,\Omega_X)$
as a complex vector space, there exists a K\"ahler form $\omega$ on $X$
such that $\tilde v([\omega])\ne0$. Consider the corresponding hypersurface
$S_{[\omega]}\subset \Def(X)$ of deformations of $X$ where $[\omega]$ stays
of type $(1,1)$ (cf.\ \ref{hypersurfaces}).
The tangent space of $S_{[\omega]}$ at
$0$ is the kernel of the map $H^1(X,\kt_X)\to H^2(X,\ko_X)\cong\IC$
induced by the product with $[\omega]\in H^1(X,\Omega_X)$.
Since $\tilde v([\omega])\ne0$, the tangent spaces $T_0S_{[\omega]}$ and
$T_0S$ have zero intersection. Hence, $S_{[\omega]}$ and
$S$ meet transversally in $0\in \Def(X)$. Shrinking $S$ we can also assume
that $S_{[\omega]}\cap S=\{0\}$.

Next, pick classes $\alpha_i\in H^2(X,\IQ)$
converging to $[\omega]$ and consider the associated hypersurfaces
$S_{\alpha_i}\subset \Def(X)$. Then the hypersurfaces $S_{\alpha_i}$ converge 
to $S_{[\omega]}$ and therefore $S_{\alpha_i}\cap S\ne\emptyset$ for $i\gg0$.
Moreover, if we choose $\alpha_i$ such that they are not of type $(1,1)$
on $X$, then $0\not\in S_{\alpha_i}\cap S$.
Hence, there exist points $t_i\in (S_{\alpha_i}\cap S)\setminus\{0\}$
converging to $0$. Using the isomorphism $H^2(\kx_t,\IZ)\cong H^2(X,\IZ)$,
the classes $\alpha_i$ are considered as rational classes of type $(1,1)$ on
$\kx_{t_i}$. Intuitively, the classes $\alpha_i$ on $\kx_{t_i}$
converge to the K\"ahler class $[\omega]$ on $X=\kx_0$ and thus should be
K\"ahler  for $i\gg0$. This can be made rigorous as follows: Fix a
diffeomorphism $\kx\cong X\times S$ compatible with the projections
to $S$. By a result of Kodaira and Spencer (Thm.\ 3.1 in \cite{KSI},
Thm.\ 15 in \cite{KSIII}) there exists a real two-form $\tilde \omega$
on $X\times S$ such that the restriction $\tilde \omega_t$ of $\tilde\omega$
to $\kx_t=X\times\{t\}$ is a K\"ahler form for all $t\in S$ and
$\tilde \omega_0=\omega$. One also finds two-forms $\tilde\omega_i$ on
$X\times S$ such that $(\tilde\omega_i)_t$ is harmonic with respect to
$\tilde\omega_t$ and $[(\tilde\omega_i)_t]=\alpha_i$ (Sect.\ 2 in
\cite{KSI}).
In particular, $(\tilde\omega_i)_{t_i}$ is a harmonic $(1,1)$-form
on $\kx_{t_i}$ representing $\alpha_i$. Since $(\tilde\omega_i)_{t_i}$
converges
to $\omega$, the harmonic $(1,1)$-form $(\tilde\omega_i)_{t_i}$ is a
K\"ahler form on $\kx_{t_i}$ for $i\gg0$. Hence, $\alpha_i$ is a Hodge class
on $\kx_{t_i}$ for $i\gg0$. In particular, $\kx_{t_i}$ is projective
for $i\gg0$.\qed

\bigskip
For the proof of the main result of this section we need to recall some
facts about positive forms and currents. Let $X$ be an irreducible symplectic
manifold and let $\omega$ be a fixed K\"ahler form on $X$.
By $\ka^{p,p}(X)$ (resp. $\ka^{p,p}(X)_\IR$) we denote the space
of smooth (real) $(p,p)$-forms on $X$. A form
$\varphi\in\ka^{p,p}(X)_\IR$ is called {\it positive} if locally it can
be written in the form
$i\alpha_1\wedge\bar\alpha_1\wedge\ldots\wedge i\alpha_p\wedge\bar\alpha_p$,
where the $\alpha_i$'s are smooth $(1,0)$-forms. For $p=2n-1$ we introduce
the convex cone $C_{pos}\subset\ka^{2n-1,2n-1}(X)_\IR$ that is spanned
by the positive forms.

\begin{remark} --- Of course, $\omega$ itself is a positive $(1,1)$-form.
Conversely, if $\varphi$ is a $d$-closed positive $(1,1)$-form, then
$\varphi+\varepsilon\omega$ is a K\"ahler form for any $\varepsilon>0$.
\end{remark} 

For the next lemma recall that $\sigma^{n-1}$ defines an isomorphism
of complex bundles ${T^*_X}^{1,0}\cong {T^*_X}^{2n-1,0}$ (use the argument
in \ref{irredcoh}).
Analogously, $\bar\sigma^{n-1}$ defines an isomorphism
${T^*_X}^{0,1}\cong {T^*_X}^{0,2n-1}$.
Passing to form valued sections of these bundles one obtains
isomorphisms ($p\geq0$):
$$\sigma^{n-1}:\ka^{1,p}(X)\cong\ka^{2n-1,p}(X)~~~~~{\rm and}~~~~~
\bar\sigma^{n-1}:\ka^{p,1}(X)\cong\ka^{p,2n-1}(X)$$ and, since
$\sigma\bar\sigma$ is real,
also $(\sigma\bar\sigma)^{n-1}:\ka^{1,1}(X)_\IR\cong \ka^{2n-1,2n-1}(X)_\IR.$

\begin{lemma}\label{posi} --- Let $\psi\in \ka^{1,1}(X)_\IR$. If
$(\sigma\bar\sigma)^{n-1}\psi\in\ka^{2n-1,2n-1}(X)_\IR$ is positive,
then also $\psi$ is positive.
\end{lemma}

\prf If $(\sigma\bar\sigma)^{n-1}\psi$ is positive it can
locally be written as
$$\begin{array}{rcl}
(\sigma\bar\sigma)^{n-1}\psi&=&i\alpha_1\wedge\bar\alpha_1\wedge\ldots\wedge i\alpha_{2n-1}\wedge\bar\alpha_{2n-1}\\
&=&i^{2n-1}(-1)^{(2n-1)(n-1)}\alpha_1\wedge\ldots\wedge\alpha_{2n-1}\wedge\bar\alpha_1\wedge\ldots\wedge\bar\alpha_{2n-1}.\\
\end{array}$$
First, observe that $i^{2n-1}(-1)^{(2n-1)(n-1)}=i$. Second, there exists
a $(1,0)$-form $\beta$ such that $\sigma^{n-1}\wedge\beta=\alpha_1\wedge\ldots
\wedge\alpha_{2n-1}$. Hence, $(\sigma\bar\sigma)^{n-1}\psi=i(\sigma\bar
\sigma)^{n-1}(\beta\wedge\bar\beta)$. Therefore,
$\psi=i\beta\wedge\bar\beta$, i.e.\ $\psi$ is positive.\qed

\bigskip
The choice of the K\"ahler form $\omega$ induces a Hodge decomposition
$$\ka^{p,p}(X)=\kh^{p,p}(X)\oplus(\im (d)\oplus \im (d^*))\cap\ka^{p,p}(X).$$
The space $\kh^{p,p}(X)$ of harmonic $(p,p)$-forms is naturally isomorphic
to $H^{p,p}(X)$. Since $\sigma\bar\sigma=\omega_J^2+\omega_K^2$ is harmonic,
where $(I,J,K)$ is the hyperk\"ahler structure associated with $(\omega,I)$
(cf.\ \ref{Comparison}), the isomorphism in \ref{irredcoh}
can be understood as the isomorphism given by
$$\begin{array}{ccc}
\kh^{1,1}(X)_\IR&\cong&\kh^{2n-1,2n-1}(X)_\IR\\
\alpha&\mapsto&(\sigma\bar\sigma)^{n-1}\alpha\\
\end{array}$$

By definition a $(1,1)$-{\it current} $T$ is a continuous linear map
$\ka^{2n-1,2n-1}(X)\to\IC$. The topology on $\ka^{2n-1,2n-1}(X)$
is the one induced by the Hodge metric.
$T$ is {\it positive} if $T(\varphi)\geq0$ for
all $\varphi\in C_{pos}$. The current $T$ is {\it closed} if $T(d\gamma)=0$ for all
$\gamma$ or, equivalently, if $T$ factorizes over
$\ka^{2n-1,2n-1}(X)/(\im (d)\cap\ka^{2n-1,2n-1}(X))$.
Any closed $(1,1)$-current
$T$ gives rise to a cohomology class $[T]\in H^{1,1}(X)$ which is given by
restricting $T$ to $\kh^{2n-1,2n-1}(X)$ and identifying the dual space
of $\kh^{2n-1,2n-1}(X)$ with $\kh^{1,1}(X)$. A current $T$ is {\it real}
if it is the natural extension of a continuous linear map
$\ka^{2n-1,2n-1}(X)_\IR\to\IR$.

\begin{proposition} --- Let $X$ be an irreducible symplectic manifold and let
$\alpha\in H^{1,1}(X)_\IR$ such that $q_X(\alpha,\,.\,)$ is positive on
the K\"ahler cone $\kk_X$. Then there exists a closed positive real
$(1,1)$-current $T$ such that $[T]=\alpha$.
\end{proposition}

\prf The proof is inspired by an argument of Peternell in \cite{Peternell}.

Let $A:=\ka^{2n-1,2n-1}(X)_\IR/(\im (d)\cap\ka^{2n-1,2n-1}(X)_\IR)$,
let $C$ be the image of $C_{pos}$ under the projection
$\ka^{2n-1,2n-1}(X)_\IR\to A$, and let $B:=\kh^{2n-1,2n-1}(X)_\IR$.
Using Hodge decomposition one has
$A=\kh^{2n-1,2n-1}(X)_\IR\oplus (\im (d^*)\cap\ka^{2n-1,2n-1}(X)_\IR)$.
Thus $B$ can be considered as a (finite dimensional)
subspace of $A$. The class $\alpha$ defines a continuous linear
map  $f_0:B\to \IR$. A closed positive real $(1,1)$-current $T$ with
$[T]=\alpha$ is a continuous linear extension
$f:A\to \IR$ of $f_0$ with $f|_C\geq0$.

The existence of such an extension is due to the following general
fact: If $A$ is a topological vector space, $B\subset A$ a subspace, and
$C\subset A$ a convex cone such that $C\cap(-C)=0$ and
$B\cap \stackrel{{\scriptscriptstyle o}}{C} \ne \emptyset$, then any continuous
linear function $f_0:B\to \IR$ with $f_0|_{B\cap C}\geq0$
can be extended to a continuous linear function
$f:A\to \IR$ with $f|_C\geq0$. (\cite{Bourbaki}, Ch.II, Sect.\ 3).

Let us now verify the assumptions in our situation: Clearly,
$C$ is a convex cone as it is the image of $C_{pos}$.
If $\varphi,-\varphi'\in C_{pos}$ with $\varphi=\varphi'+d\gamma$ then the positivity
of $\omega$ implies $0\leq\int\varphi\omega=\int\varphi'\omega\leq0$
and hence $\int\varphi\omega=\int\varphi'\omega=0$.
Since $\omega$ is K\"ahler,
we get $\varphi=\varphi'=0$. Therefore, $C\cap(-C)=0$.
The intersection $B\cap \stackrel{\scriptscriptstyle o}{C}$ is not empty as it contains the class
$\omega^{2n-1}$ which is an inner point of $C_{pos}$ (cf.\ \cite{Harvey}).
Last but not least, one checks $f_0|_{B\cap C}\geq0$.
Indeed, if $\varphi\in \kh^{2n-1,2n-1}(X)_\IR$ is positive,
then there exists $\psi\in\kh^{1,1}(X)_\IR$ with
$(\sigma\bar\sigma)^{n-1}\psi=\varphi$ and $\psi $ is positive by
Lemma \ref{posi}. Hence, $\psi+\varepsilon\omega$ is K\"ahler
for all $\varepsilon>0$. The assumption on $\alpha$ now
implies $q_X(\alpha,[\psi+\varepsilon\omega])>0$ for all $\varepsilon>0$
and hence $q_X(\alpha,[\psi])\geq0$. Therefore,
$\int\varphi\alpha=\int(\sigma\bar\sigma)^{n-1}\psi\alpha=
cq_X([\psi],\alpha)\geq0$ for some positive
constant $c$. Hence, $f_0(\varphi)\geq0$.\qed 

\begin{remark} --- {\it i)} Demailly \cite{Demailly} introduced the notion
of pseudo-effective classes: A class $\alpha\in H^{1,1}(X)$ is called
pseudo-effective if it can be represented by a closed
positive $(1,1)$-current. They span a closed cone $H^{1,1}_{psef}(X)\subset
H^{1,1}(X)$. The proposition says that the cone dual to the K\"ahler cone
is contained in $H^{1,1}_{psef}(X)$.

{\it ii)} By \ref{posconedef} any class $\alpha\in \kc_X$ is positive on $\kk_X$.
Hence $\kc_X\subset H^{1,1}_{psef}(X)$. In fact, as $\kc_X$ is an open
cone, any class $\alpha\in\kc_X$ can be represented by a closed positive real
$(1,1)$-current $T$ such that $T-\varepsilon\omega$ is still positive
for some $\varepsilon>0$. 
\end{remark}

As an immediate consequence of \cite{Demailly} we obtain

\begin{corollary}\label{demcor} --- If $X$ is a projective
irreducible symplectic manifold and $L$ is a line bundle
such that $c_1(L)\in\kc_X$, then $L$ has maximal Kodaira-dimension, i.e.\
$h^0(L^k)$ grows like $c\cdot k^{2n}$ for some $c>0$.
In particular, $\kc_X\cap H^2(X,\IZ)$ is contained in the effective cone. \qed
\end{corollary}

The main result of this section is the following 

\begin{theorem}\label{proj}---
Let $X$ be an irreducible symplectic manifold. Then
$X$ is projective if and only if there exists a line bundle $L$ on $X$
with $q_X(c_1(L))>0$.
\end{theorem}

\prf If $X$ is projective then $c_1(L)$ of an ample line bundle $L$
is a K\"ahler class and, therefore, $q_X(c_1(L))>0$
(cf.\ \ref{quadraticform}).

Conversely, if $X$ admits a line bundle $L$ with $c_1(L)\in\kc_X$, then
there exists a closed $(1,1)$-current $T$ representing $c_1(L)$ such that
$T-\varepsilon\omega$ is positive for some $\varepsilon>0$.
By a result of Ji-Shiffman \cite{Shiffman} and Bonavero \cite{Bonavero}
this implies that $X$ is Moishezon and, since it is also
K\"ahler, that it is projective.\qed
 
\begin{remark} ---
{\it i)} The theorem only asserts that $X$ is projective, but not that
any line bundle $L$ with positive square $q_X(c_1(L))$ is ample
(this is not even true for surfaces).
A criterion for the ampleness of a line bundle will be discussed
in Sect.\ \ref{amplecone}. However, if $\Pic(X)\cong L\cdot \IZ$
with $q_X(c_1(L))>0$, then either $L$ or $L\dual$ is ample. 

{\it ii)} In the original approach to prove the theorem
I tried to avoid the result of
Bonavero and Ji-Shiffman, which in turn relies on Demailly's very complicated
holomorphic Morse inequalities. The argument should
use the relatively
easy Corollary \ref{demcor} on projective deformations and some kind
of semi-continuity argument.
But there are still some technical problems. Also note that
if one is willing to use the holomorphic Morse inequalities (or rather
their singular version) then one can in fact see that Corollary
\ref{demcor} works without the projectivity assumption.

{\it iii)} The techniques above can also be used to prove a result for K\"ahler
surfaces I was not aware of before
and which I could not find in the literature.
However, I believe that going through the classification of surfaces an easier
proof, not using singular Morse inequalities, should exist. The result
is: A compact K\"ahler surface is projective if and only if the dual of the
K\"ahler cone, i.e.\ the elements with
positive intersection with any K\"ahler class,
contains an inner integral point.
\end{remark}

\section{Birational Manifolds}\label{biratman}

This section deals with the relation between birational irreducible symplectic
manifolds and non-separated points in the moduli space $\gM_\Gamma$
(cf.\ \ref{moduli}). Already in dimension two, i.e.\ for
K3 surfaces, the moduli space of marked irreducible symplectic manifolds
is not separated (i.e.\ non-Hausdorff). But there, two non-separated points
always correspond to just one K3 surface with two different markings. The
situation is more subtle in higher dimensions: One easily
generalizes the `Main Lemma' of Burns and Rapoport \cite{BurnsRapoport},
the algebraic version of which is \cite{MM},
to the effect that the underlying manifolds $X$ and $X'$ of two non-separated
points $(X,\varphi),(X',\varphi')\in\gM_\Gamma$ are birational (Theorem
\ref{BRMain}).
But contrary to the two-dimensional case, this does not imply
that $X$ and $X'$ are isomorphic. In fact, it is a standard procedure
to construct via certain birational transformations
out of one irreducible symplectic manifold new ones (cf.\ \ref{exbybir}).
It is the
goal of this section to prove that in general two birational irreducible
symplectic manifolds correspond to non-separated points in $\gM_\Gamma$
(Theorem \ref{birat}). As the method is completely algebraic the result
is limited to the case of projective manifolds (but see Sect.\
\ref{remarks} for the
relation between a conjectural Global Torelli Theorem and this
result for non-projective manifolds). This part is a continuation
of the predecessor \cite{Huybrechts} of this paper, where the result was
proved under an additional assumption on the codimension of the
exceptional locus. For some details of the proof we will refer to
\cite{Huybrechts}.

The result has two applications. The first, to be considered in this
section, concerns the classification of known examples of irreducible
symplectic manifolds (cf.\ Sect.\ \ref{examples}). Corollary \ref{onlytwo}
shows that all known examples of irreducible symplectic manifolds are
deformation equivalent to one of the two standard examples: to the Hilbert
scheme of points on a K3 surface or to a generalized Kummer variety (see
Sect.\ \ref{examples} for their definition). In particular, in any real
dimension $4n>4$ we know exactly two compact differentiable manifolds
admitting an irreducible hyperk\"ahler metric. (For $n=1$
one knows that there exists exactly one: the real manifold underlying
a K3 surface.) In the light of this result the list of
examples of irreducible symplectic manifolds seems rather short.
It might be noteworthy
that the general results of this paper, e.g.\ the description of the K\"ahler
cone \ref{kaehlerconethm},
the projectivity criterion \ref{proj} or the surjectivity of the period map \ref{surjper},
are by no means trivial (or a direct consequence of the K3 surface theory)
for the deformations of the two standard
series $\Hilb^n(S)$ and ${\rm K}^{n+1}(A)$.
In fact, I do not see how one could
possibly simplify the proofs for these special cases.

The second application concerns a conjectural
`Global Torelli Theorem' for higher
dimensional irreducible symplectic manifolds.
As this problem is intimately related to questions about 
the period map we postpone the discussion until Sect.\ \ref{remarks}.

Before going into the subject we say a few words about what `non-separated'
means in practice. To this end we formulate the following
(almost tautological):

\begin{lemma}\label{lemsep}---
Let $X$ and $X'$ be irreducible symplectic manifolds, let $\kx\to S$
and $\kx'\to S$ be deformations of $X$ and $X'$, respectively, and
let $V\subset S$ be an open (in the analytic topology) non-empty subset
such that 
{\it i)} $S$ is one-dimensional and {\it ii)}
$0\in\partial V$, and {\it iii)} $\kx|_V\cong
\kx'|_V$ (compatible with the projections to $S$). Then
there exist markings $\varphi$ and $\varphi'$ of $X$ and $X'$, respectively,
such that $(X,\varphi),(X',\varphi')\in\gM_\Gamma$ are non-separated.
Conversely, if $(X,\varphi),(X'\varphi')\in\gM_\Gamma$ are two non-separated
points in
$\gM_\Gamma$, then there exist deformations $\kx\to S$ and $\kx'\to S$
and $V\subset S$ satisfying {\it i)}, {\it ii)}, and {\it iii)}
such that for $t\in V$ the natural isomorphisms
$H^2(X,\IZ)\cong H^2(\kx_t,\IZ)\cong H^2(\kx'_t,\IZ)\cong H^2(X',\IZ)$ are
compatible with $\varphi$ and $\varphi'$.\qed
\end{lemma}

\begin{remark}--- The arguments in the proof of Proposition \ref{Weylprop}
show that in part two of the lemma one can in fact arrange things such
that $V=S\setminus\{0\}$. See Remark \ref{topoo}.
\end{remark}

Let us now come to a straightforward generalization of the `Main Lemma'
in \cite{BurnsRapoport}. The proof follows closely Beauville's expos\'e
in \cite{Periodes}.

\begin{theorem}\label{BRMain}---
If $(X,\varphi), (X',\varphi')\in\gM_\Gamma$ are non-separated points
in the moduli space of marked irreducible symplectic manifolds, then $X$
and $X'$ are birational.
\end{theorem}

\prf Consider deformations $\kx\to S$ and $\kx'\to S$ of $\kx_0=X$ and
$\kx'_0=X'$, respectively, and $V\subset S$
as in \ref{lemsep}. In particular, we have
for $t\in V$ the canonical isomorphism
$H^2(X,\IZ)\cong H^2(\kx_t,\IZ)\cong H^2(\kx'_t,\IZ)\cong H^2(X',\IZ)$
which is compatible with $\varphi$ and $\varphi'$. Pick a sequence
$t_i\in V$ converging to $0\in S$ and consider the corresponding
isomorphisms $f_i:\kx_{t_i}\cong \kx'_{t_i}$ and their graphs
$\Gamma_i\subset \kx_{t_i}\times \kx'_{t_i}$. By \cite{KSI} there exist
K\"ahler forms $\omega_t$ and $\omega'_t$ on the fibres $\kx_t$,
respectively $\kx'_t$, for all $t\in S$ depending continuously (this is enough)
on $t$. The volume of $\Gamma_{t_i}$ with respect
to these K\"ahler forms, can be computed by

$$\begin{array}{ccl}
vol(\Gamma_i)&=&\displaystyle{\int_{\kx_{t_i}}(\omega_{t_i}+f_i^*\omega'_{t_i})^{2n}}\\
&=&\displaystyle{\int_{\kx_{t_i}}([\omega_{t_i}]+f_i^*[\omega'_{t_i}])^{2n}}\\
&=&\displaystyle{\int_X([\omega_{t_i}]+({\varphi'}^{-1}\circ\varphi)[\omega'_{t_i}])^{2n}}.\\
\end{array}$$

(Use $f_i^*=({\varphi'}^{-1}\circ \varphi)$ via the isomorphisms
$H^2(X,\IZ)\cong H^2(\kx_t,\IZ)$ and $H^2(X',\IZ)\cong H^2(\kx'_t,\IZ)$.)

Hence $vol(\Gamma_i)$ converges to
$\int_X([\omega_{0}]+({\varphi'}^{-1}\circ\varphi)[\omega'_{0}])^{2n}<\infty$.
By a result of Bishop \cite {Bishop}
the boundedness of the volume of $\Gamma_{i}$ is enough to conclude the
existence of a limit cycle $\Gamma_\infty\subset X\times X'$ with the same
cohomological
properties as the $\Gamma_i$'s. In particular,
$[\Gamma_\infty]\in H^{4n}(X\times X',\IZ)$ satisfies
$p_*[\Gamma_\infty]=[X]\in H^0(X,\IZ)$;
$p'_*[\Gamma_\infty]=[X']\in H^0(X',\IZ)$; and
$p'_*([\Gamma_\infty].p^*\alpha)=({\varphi'}^{-1}\circ\varphi)(\alpha)$
for all $\alpha\in H^2(X,\IZ)$.
Here, $p$ and $p'$ denote the two projections from $X\times X'$.
Splitting $\Gamma_\infty$ into its irreducible components and using the first
two properties we have either
\begin{itemize}
\item[$\bullet$] $\Gamma_\infty=Z+\sum Y_i$, where $p:Z\to X$ and
$p':Z\to X'$ are generically one-to-one, or
\item[$\bullet$]  $\Gamma_\infty=Z+Z'+\sum Y_i$, where $p:Z\to X$ and
$p':Z'\to X'$ are generically one-to-one, but neither
$p':Z\to X'$ nor $p:Z'\to X$ is generically finite.
\end{itemize}
In both cases, $p_*[Y_i]=0$ and $p'_*[Y_i]=0$.

The second possibility can be excluded: If $\sigma$ and $\sigma'$
are non-trivial holomorphic two-forms on $X$ and $X'$ respectively, then
$({\varphi'}^{-1}\circ\varphi)([\sigma])=\lambda[\sigma']$ for some
$\lambda\ne0$. Indeed, since the period domain, i.e.\ the quadric
$Q\subset \IP(\Gamma_\IC)$, is separated, ${\varphi'}^{-1}\circ\varphi$ is
compatible with the natural Hodge structures on $H^2(X,\IZ)$ and
$H^2(X',\IZ)$. If
$\Gamma_\infty=Z+Z'+\sum Y_i$, then 

$$\begin{array}{ccl}
0&\ne&\displaystyle{\lambda\int_{X'}(\sigma'\bar\sigma')^n=\int_{X'}({\varphi'}^{-1}\circ\varphi)[\sigma]({\sigma'}^{n-1}{\bar\sigma'}{}^{n})}\\
&=&\displaystyle{\int_{X'}p'_*([\Gamma_\infty].p^*[\sigma])({\sigma'}^{n-1}{\bar\sigma'}{}^n)=\int_{X\times X'}[\Gamma_\infty].p^*[\sigma].{p'}^*[{\sigma'}^{n-1}{\bar\sigma'}{}^n]}\\
&=&\displaystyle{\int_Zp^*[\sigma].{p'}^*[{\sigma'}^{n-1}{\bar\sigma'}{}^n]+\int_{Z'}p^*[\sigma].{p'}^*[{\sigma'}{}^{n-1}{\bar\sigma'}{}^n]+\sum\int_{Y_i}p^*[\sigma].{p'}^*[{\sigma'}{}^{n-1}{\bar\sigma'}{}^n]}.\\
\end{array}$$

The first and third term vanish, because $p'(Z)$ and $p'(Y_i)$ are
of dimension $<2n$, but ${\bar\sigma'}{}^n$ is a $(0,2n)$-form on $X'$.
For simplicity assume $Z'$ smooth (otherwise pass to a desingularization).
Since $Z'\to X'$ is birational, $H^0(Z',\Omega^2_{Z'})={p'}^*({\sigma'})|_{Z'}\cdot\IC$. On the other hand, $p^*(\sigma)|_{Z'}\in H^0(Z',\Omega_{Z'}^2)$.
Since $p^*(\sigma)|_{Z'}$ is everywhere degenerate, but ${p'}^*(\sigma')|_{Z'}$
is at least generically non-degenerate, we must have $p^*(\sigma)|_{Z'}=0$.
Thus, also the second term vanishes. Contradiction.

Therefore, only the decomposition $\Gamma_\infty=Z+\sum Y_i$ can occur. Since
$Z\to X$ and $Z\to X'$ are generically one-to-one, $X$ and $X'$ are
birational.\qed

\bigskip
Note that the birational correspondence between $X$ and $X'$ constructed
this way does not, in general, induce ${\varphi'}^{-1}\circ \varphi$ on $H^2$.
The theorem (or rather its proof) has various interesting consequences
and, with the proof still in mind, the reader may wish to have a look
at Corollary \ref{corBRMain} already at this point in the discussion.

\bigskip

Let us now come to the converse: Are two birational irreducible symplectic
manifolds non-separated in their moduli space?
Recall that for an appropriate choice of the markings
the periods of two birational irreducible symplectic manifolds are
equal (cf.\ \ref{hodgeunderbir}).
Therefore, their period points in $\IP(\Gamma_\IC)$ coincide.

Let us consider the following situation:

\bigskip
\refstepcounter{theorem}\label{assump}{\bf \thetheorem} ---
{\it i)}
$X$ and $X'$ are irreducible symplectic manifolds.
{\it ii)}
$X'$ is projective.
{\it iii)}
There exists a birational map $f:X - - \to X'$.

\bigskip
Let $U\subset X$ and $U'\subset X'$ be the maximal open subsets where $f$,
respectively $f^{-1}$, are regular. Then $U\cong U'$ and $\codim(X\setminus U),
\codim(X'\setminus U')\geq2$. This is a general fact for varieties
with nef canonical divisor, but see \cite{Huybrechts}.
Since $X$ is K\"ahler and birational to the projective manifold
$X'$, also $X$ is projective (cf.\ \ref{Moish}). For the algebraic Picard
groups we have canonical isomorphisms
$\Pic(X)\cong\Pic(U)\cong\Pic(U')\cong\Pic(X')$. If $L\in\Pic(X)$ we denote
by $L'$ the associated line bundle on $X'$, i.e.\ $L'$ is
the line bundle on $X'$ such that
$L|_U\cong L'|_{U'}$. In the sequel, $L'$ will usually be ample. 
If $L$ is ample as well, then the isomorphism $U\cong U'$ can be extended to an
isomorphism of $X$ and $X'$ and there is nothing to show.

The following proposition was proved in \cite{Huybrechts}.

\begin{proposition}\label{biratdefo}--- Under the assumptions
\ref{assump}, let $L'$ be ample and
let $\pi:(\kx,\kl)\to S$ be a deformation of
$(\kx_0,\kl_0)=(X,L)$ over a smooth one-dimensional base $S$.
If for any $n\gg0$ the dimension $h^0(\kx_t,\kl_t^n)$ is constant in a
neighbourhood of $0\in S$, then, after shrinking $S$ if necessary,
there exists a deformation $(\kx',\kl')\to S$ of $(\kx'_0,\kl'_0)=(X',L')$
and an $S$-birational correspondence between $\kx$ and $\kx'$ respecting
$\kl$ and $\kl'$.\qed
\end{proposition}

The idea of the proof is as follows: The linear system $|\kl^m_t|$
($m\gg0$)
defines a rational map from $\kx$ to $\IP_S((\pi_*\kl^m)\dual)$
and we define $\kx'$ as the closure of the image of this rational map.
Then one shows that $\kx'$ is $S$-flat with special fibre $X'$.
The line bundle ${\kl'}^m$ is the restriction of the relative $\ko(1)$
on $\IP_S((\pi_*\kl^m)\dual)$. For details we refer to \cite{Huybrechts}.
Note that
the neighbourhood where $h^0(\kx_t,\kl_t^m)\equiv const$
might very well depend on $n$.  

The main technical problem that limited the results in \cite{Huybrechts}
to the case that $\codim(X\setminus U),\codim(X'\setminus U')\geq3$
was the assumption $h^0(\kx_t,\kl_t^n)\equiv const$, which
is easy to establish under this additional assumption on the codimension.
The projectivity criterion Theorem \ref{proj}, which was not
yet available in \cite{Huybrechts}, can be used to arrange things such that
the assumption on $h^0(\kx_t,\kl_t)$ holds.

\begin{theorem}\label{birat}---
Let $X$ and $X'$ be birational projective irreducible symplectic
manifolds. Then there exist deformations $\kx\to S$ and $\kx'\to S$ of
$\kx_0=X$ and $\kx'_0=X'$, respectively, such that
\begin{itemize}
\item[$\bullet$] $S$ is smooth and one-dimensional.
\item[$\bullet$] There exists an $S$-isomorphism $\kx|_{S\setminus\{0\}}\cong
\kx'|_{S\setminus\{0\}}$.
\end{itemize}
\end{theorem}

\prf By \ref{Todd}
there exist constants $(a_i)$ and $(a'_i)$ ($i=1,\ldots,n$) such
that for any line bundle $L$ on $X$ (resp.\ $L'$ on $X'$) the
Hirzebruch-Riemann-Roch formula can be written as
$$\chi(X,L)=\sum_{i=1}^n\frac{a_i}{(2i)!}q_X(c_1(L))^i\phantom{XX}
{\rm resp.}\phantom{XX}\chi(X',{L'})=\sum_{i=1}^n\frac{a'_i}{(2i)!}q_{X'}(c_1(L'))^i.$$
Without loss of generality we may assume that in the
lexicographic order $(a_i)\geq(a'_i)$, i.e.\ $\sum a_i t^i\geq\sum a_i't^i$
for $t\gg0$.
If now $L'\in\Pic(X')$ is ample and $L$ is the corresponding line bundle on
$X$, then $q_X(c_1(L))=q_{X'}(c_1(L'))>0$.
For the equality see \ref{hodgeunderbir} and for the inequality
\ref{quadraticform}.
Hence, $\chi(X,L^m)\geq\chi(X',{L'}^m)$ for $m\gg0$.

Let $(\kx,\kl)\to S$ be a deformation of $(X,L)$ over a smooth and
one-dimensional base $S$
such that $\rho(\kx_t)=1$ for general $t\in S$, e.g.\ take a
curve $S\subset \Def(X,L)$ not contained in any $\Def(X,M)$ where
$M$ is a line bundle linearly independent of $L$ (cf.\ \ref{deformationlb},
\ref{hypersurfaces}).
By Theorem \ref{proj} the positivity of
$q_{\kx_t}(c_1(\kl_t))=q_X(c_1(L))$ implies that all fibres
$\kx_t$ are projective. Hence, for general $t\in S$ either $\kl_t$ or
its dual is ample. The non-vanishing of $H^0(X,L^m)=H^0(X',{L'}^m)$ for
$m\gg0$ and the semi-continuity of $h^0(\kx_t,\kl_t)$ excludes the latter
possibility for $t$ close to
$0\in S$. Thus, we can apply
the Kodaira Vanishing Theorem to $\kl_t$:
For general $t\in S$ one has $h^0(\kx_t,\kl_t^m)=\chi(\kx_t,\kl_t^m)
=\chi(X,L^m)$ for all $m>0$.
Using semi-continuity this yields $h^0(X,L^m)\geq h^0(\kx_t,\kl_t^m)=\chi(X,L^m)$ for $m>0$.
On the other hand, $h^0(X',{L'}^m)=h^0(X,L^m)$
(use $\codim(X\setminus U),\codim(X'\setminus U')\geq2$) and
$\chi(X',{L'}^m)=h^0(X',{L'}^m)$ (use $L'$ ample on $X'$). Therefore,
$$\chi(X',{L'}^m)=h^0(X,L^m)\geq h^0(\kx_t,\kl_t^m)=\chi(X,L^m).$$
Under the assumptions $(a_i)\geq(a'_i)$ and $q_X(c_1(L))=
q_{X'}(c_1(L'))>0$ this
implies $h^0(X,L^m)=h^0(\kx_t,\kl_t^m)$ for $m\gg0$ .
Hence, for $m\gg0$ the dimension $h^0(\kx_t,\kl_t^m)$ is constant
in a neighbourhood of $0\in S$.

Therefore, Proposition \ref{biratdefo} can be applied and we find a deformation
$(\kx',\kl')\to S$ of $(X',L')$ which is $S$-birational to $\kx$. The rest
of the proof is as the one of the corresponding
theorem in \cite{Huybrechts}: For $t\in S$ general and close to $0\in S$
the birational correspondence $\kx_t- - \to \kx'_t$ is a birational
correspondence between two projective manifolds with Picard number $\rho(\kx_t)=\rho(\kx_t')=1$, which must be an isomorphism.
Hence, $\kx- - \to \kx'$ is an isomorphism on the
general fibre. Restricting to an open neighbourhood of $0\in S$ we can in
fact achieve that $\kx|_{S\setminus\{0\}}\cong\kx'|_{S\setminus\{0\}}$ (cf.\
\cite{Huybrechts}).\qed

\bigskip
Note that the theorem is known as well for non-projective irreducible
symplectic manifolds if the birational correspondence is
described by an elementary transformation \cite{Huybrechts}.

More in the spirit of the Main Lemma of Burns and Rapoport (see \ref{BRMain})
Theorem \ref{birat} can equivalently be formulated as

\bigskip
{\parindent0mm
{\bf Theorem \ref{birat}'}} ---
{\it Let $X$ and $X'$ be birational projective
irreducible symplectic manifolds. Then there exist
two markings $\varphi:H^2(X,\IZ)\cong \Gamma$ and
$\varphi':H^2(X',\IZ)\cong\Gamma$ such that
$(X,\varphi),(X',\varphi')\in\gM_\Gamma$ are non-separated points.\qed}
\bigskip

The fact that $X$ and $X'$ can be realized as the special fibres of the same
family has strong consequences, which are not at all obvious
just from the fact that they are birational.
We only mention:

\begin{corollary}---
If $X$ and $X'$ are birational projective irreducible symplectic
manifolds, then:
\begin{itemize}
\item[{\it i)}] $X$ and $X'$ are diffeomorphic.
\item[{\it ii)}] For all $k$ the weight-$k$ Hodge structures of $X$ and $X'$
are isomorphic.\qed
\end{itemize}
\end{corollary}

It is interesting to compare {\it ii)} in the corollary with
a recent result of Batyrev and Kontsevich. They show that
{\it ii)} holds true for all birational smooth projective manifolds
with trivial canonical bundle; in particular for irreducible symplectic
but also for Calabi-Yau manifolds.
Assertion {\it i)} is not expected to hold for the more general class
of projective manifolds with trivial
canonical bundle.

\bigskip

{\bf Applications}

\bigskip
In Sect.\ \ref{examples} we provided a list of the known examples
of irreducible symplectic manifolds. In most of the cases the verification
of the defining properties was reduced to either of the two standard
examples $\Hilb^n(S)$ or ${\rm K}^{n+1}(A)$, where $S$ is a K3 surface and $A$ is a
torus. In fact, modulo birational correspondence all examples
in Sect.\ \ref{examples} are deformation
equivalent to one of these two. Theorem \ref{birat} in particular
says that also birational irreducible symplectic manifolds are
deformation equivalent. In particular, we wish to mention

\begin{corollary}--- If $S$ is a K3 surface and $v=(v_0,v_1,v_2)$ is a Mukai
vector with $v_1$ primitive, then for any $v$-generic polarization $H$
the moduli space $M_H(v)$ of semistable sheaves with Mukai vector $v$ is an
irreducible symplectic manifold that is
deformation equivalent to $\Hilb^n(S)$, where
$n=(v,v)+2$. (For the notation see Sect.\ \ref{examples}, \cite{OG} or
\cite{HL}.)\qed
\end{corollary}

More generally, we formulate

\begin{corollary}\label{onlytwo}---
All examples of irreducible symplectic manifolds in Sect.\ \ref{examples}
are deformation equivalent (and hence diffeomorphic) either
to $\Hilb^n(S)$ or to ${\rm K}^{n+1}(A)$, where $S$ is a K3 surface and $A$
is a torus.
\end{corollary}

\prf The only problem is the projectivity assumption in Theorem \ref{birat}.
But since all known examples of birational correspondences between
non-projective
irreducible symplectic manifolds are described
in terms of elementary transformations
along projective bundles, the result follows from \cite{Huybrechts}, where
we gave an easy proof of Theorem \ref{birat} for elementary transformations
without assuming the projectivity of $X$ and $X'$.\qed


\section{An Analogue of the Weyl-Action}\label{weylaction}

The main obstacle to generalize the techniques in the theory
of K3 surfaces to
higher dimensions, besides the missing Global Torelli Theorem,
is the absence of a Weyl-group, i.e.\ the group
of automorphisms of the second cohomology generated by reflections
in hypersurfaces orthogonal to some $(-2)$-class. To a certain extent, the following proposition
(and its corollary) is a good replacement for the fact that the Weyl-group
acts transitively on the set of chambers (see also Remark \ref{Weyl}).
It has immediate consequences such as the description of
the K\"ahler cone of very general irreducible symplectic
manifolds (see \ref{Kaehlerconeforperiod} and \ref{Kaehlerconeforperiod2}).
The question will be further pursued in Sect.\ \ref{kaehlercone}.
A posteriori the not very concrete `generality' assumption in
\ref{Weylprop} can be made more specific (Theorem \ref{Delignesuggestion}).
This was pointed out to me by Deligne.
 
\begin{proposition}\label{Weylprop}---
Let $(X,\varphi)\in\gM_\Gamma$ be a marked irreducible symplectic manifold.
Assume $\alpha\in\kc_X$ is general, i.e.\ $\alpha$ is contained in the
complement of countably many nowhere dense
closed subsets. Then there exists a point
$(X',\varphi')\in\gM_\Gamma$, which cannot be separated from $(X,\varphi)$
such that $({\varphi'}^{-1}\circ\varphi)(\alpha)\in H^2(X',\IR)$
is a K\"ahler class.
\end{proposition}

Together with Theorem \ref{BRMain} (or rather its proof) the proposition shows

\begin{corollary}\label{corBRMain} ---
Let $X$ be an irreducible symplectic manifold of dimension
$2n$ and let $\alpha\in\kc_X$
be general. Then there exists another irreducible symplectic manifold
$X'$ together with an effective cycle $\Gamma:=Z+\sum Y_i\subset X\times X'$
of dimension $2n$ satisfying the following conditions:
\begin{itemize}
\item[{\it i)}] $Z$ defines a birational map $X- - \to X'$.
\item[{\it ii)}] The projections $Y_i\to X$ and $Y_i\to X'$ have positive fibre
dimension.
\item[{\it iii)}] $[\Gamma]_*:H^2(X,\IZ)\to H^2(X',\IZ)$ defines an isomorphism
of Hodge structures compatible with $q_X$ and $q_{X'}$. Moreover,
$[\Gamma]_*\circ[\Gamma]_*$ is the identity on $H^2(X,\IZ)$, resp.\
$H^2(X',\IZ)$.
\item[{\it iv)}] $[\Gamma]_*(\alpha)\in\kk_{X'}$.\qed
\end{itemize}
\end{corollary}

{\it Proof of Proposition \ref{Weylprop}.}
First,  let $\alpha\in \kc_X$ be completely arbitrary. Choose a K\"ahler class
$\gamma\in\kk_X\subset\kc_X$ and let $\beta:=(1+\varepsilon)\alpha-\varepsilon\gamma\in\kc_X$, where $0<\varepsilon \ll1$. Then, pick a sequence $\beta_i\in H^2(X,\IQ)$ converging to $\beta$ such that $\beta_i\not\in H^{1,1}(X)$
(or, equivalently, $\beta_i\not\in\kc_X$). The associated hypersurfaces $S_{\beta_i}\subset \Def(X)$ (for the notation see \ref{hypersurfaces})
converge to $S_\beta\subset \Def(X)$. Since $S_\beta\ne\emptyset$, also
$S_{\beta_i}\ne\emptyset$ for $i\gg0$. Let $U\subset \kc_X$ be an open neighbourhood of $\alpha$ such that for all $\alpha'\in U$ the class
$((1+\varepsilon)\alpha'-\beta)/\varepsilon$ is still K\"ahler.
Such an open neighbourhood exists, because $\kk_X\subset \kc_X$ is open.
For any $\alpha'\in U$ let $\kx(\alpha')\to T(\alpha')$ be the deformation
obtained by restricting the universal deformation to
$T(\alpha'):=\Def(X)\cap\kp^{-1}(\IP(\varphi(F(\alpha')_\IC)))$ (see also
\ref{twistor2}).
Since $\alpha',\beta\in\kc_X$ and, hence, $q_X(\alpha',\beta)\ne0$, the
curve $T(\alpha')$ is not contained in $S_\beta$ (cf.\ \ref{pairing},
\ref{hypersurfaces}). Therefore,
$T(\alpha')\cap S_\beta$ is zero-dimensional and, moreover,
non-empty, as the orign $0$ is contained in
$T(\alpha')\cap S_\beta$. Hence, $T(\alpha')\cap S_{\beta_i}$ is
zero-dimensional and non-empty for $i\gg0$ as well. Thus, we find points
$t_i\in T(\alpha')\cap S_{\beta_i}$ approaching $0$.

We need the following

\bigskip
{\it Claim} --- For general $\alpha'\in U$ the points $t_i\in T(\alpha')\cap
S_{\beta_i}$ are general points of $S_{\beta_i}$ in the sense that $H^{1,1}(\kx(\alpha')_{t_i})_\IQ$ is spanned by $\beta_i$ or, equivalently,
that $\rho(\kx(\alpha')_{t_i})=1$.

\bigskip
{\it Proof of the Claim.} Consider the map $\psi_i:U\to S_{\beta_i}$
that sends $\alpha'\in U$ to the intersection $\{t_i\}=T(\alpha')\cap S_{\beta_i}$, which, at least locally, consists of a single point. The map $\psi_i$
is constant along the orbits of the natural $\IR^*$-action on
$U\subset H^{1,1}(X)_\IR$. Thus $d\psi_i:T_{\alpha'}U\to T_{t_i} S_{\beta_i}$
factorizes through
$T_{\alpha'}U=H^{1,1}(X)_\IR\to H^{1,1}(X)_\IR/\alpha'\cdot\IR
\to T_{t_i} S_{\beta_i}\subset H^{1,1}(\kx(\alpha')_{t_i})$.
Using the period map, which is holomorphic, it is easy
to extend the map $\psi_i$ to a holomorphic map
$\tilde\psi_i:\tilde U\to S_{\beta_i}$, where $\tilde U$ is an open subset of
the complex vector space $H^{1,1}(X)$ with $\tilde U\cap H^{1,1}(X)_\IR=U$.
Again, $\tilde \psi_i$ is constant along the orbits of the natural
$\IC^*$-action on $\tilde U\subset H^{1,1}(X)$ and, therefore,
its tangent map factorizes as follows
$T_{\alpha'}\tilde U\cong H^{1,1}(X)\to H^{1,1}(X)/\alpha'\cdot\IC\to
T_{t_i}S_{\beta_i}\subset H^{1,1}(\kx(\alpha')_{t_i})$.

Next, one proves that the map $\tilde\psi:\tilde U\to S_{\beta_i}$ is injective
modulo the $\IC^*$-action. Indeed, if $H^{2,0}=\sigma\cdot\IC$, then
$F(\alpha')_\IC=\sigma\cdot\IC\oplus\bar\sigma\cdot\IC\oplus\alpha'\cdot\IC$.
Therefore, for linearly independent $\alpha',\alpha''\in\tilde U\subset
H^{1,1}(X)$ the two
planes $\IP(\varphi(F(\alpha')_\IC)),\IP(\varphi(F(\alpha'')_\IC))\subset
\IP(\Gamma_\IC)$ intersect in $\IP(\varphi(\sigma\cdot\IC\oplus\bar\sigma\cdot\IC))$.
The latter space is a projective line which meets the period domain $Q\subset \IP(\Gamma_\IC)$ in exactly two points, namely
$\IP(\varphi(\sigma\cdot\IC\oplus\bar\sigma\cdot\IC))\cap Q=\{\varphi[\sigma],\varphi[\bar\sigma]\}$. Shrinking $\Def(X)$ such that
$\varphi[\bar\sigma]\not\in\kp(\Def(X))$,
one has in particular
$\varphi[\bar\sigma]\not\in\kp(S_{\beta_i})$ and, therefore,
$T(\alpha')\cap S_{\beta_i}\ne T(\alpha'')\cap S_{\beta_i}$ for linearly
independent $\alpha'$ and $\alpha''$. Thus $\tilde\psi:\tilde U\to S_{\beta_i}$
is injective modulo $\IC^*$. 

But if $\tilde\psi_i:\tilde U\to S_{\beta_i}$
is injective modulo the $\IC^*$-action, then its tangent map
$d\tilde\psi_i:H^{1,1}(X)/\alpha'\cdot\IC\to T_{t_i}S_{\beta_i}$ is injective
for general $\alpha'\in \tilde U$.
In fact, the set where $d\tilde \psi_i$
fails to be injective is a complex-analytic set and, therefore,
cannot contain open parts of $U$. Thus, even for
general $\alpha'\in U$ the tangent map can be assumed to be injective.
Since both spaces $H^{1,1}(X)/\alpha'\cdot\IC$ and $T_{t_i}S_{t_i}$
have dimension $h^{1,1}(X)-1$, the tangent map $d\tilde\psi_i$ at such a point
must be bijective. In particular, $\im(d\tilde\psi_i)$ is not contained in any
$T_{t_i}S_\delta\cap T_{t_i}S_{\beta_i}$ for any $(1,1)$-class $\delta$ on
$\kx(\alpha')_{t_i}$ that is linearly independent of $\beta_i$. Since $d\tilde\psi_i$ is $\IC$-linear and $H^{1,1}(X)_\IR$ spans $H^{1,1}(X)$, this
also shows that $\im(d\psi_i)$ is not contained in any such
$T_{t_i}S_\delta\cap T_{t_i}S_{\beta_i}$ . Hence, the image of
$\psi_i:U\to S_{\beta_i}$ is not contained in any hypersurface
$S_\delta\cap S_{\beta_i}$ with $\delta$ linearly independent of $\beta_i$.
As there are only countably many $\beta_i$'s and $\delta$'s
to be considered, one can assume that for the
general $\alpha'\in U$ the intersection $T(\alpha')\cap S_{\beta_i}$ is not
contained in any $S_\delta$ for $\delta$ linearly independent of $\beta_i$.
In other words, if $\{t_i\}=T(\alpha')\cap S_{\beta_i}$ with $\alpha'$ general,
then $\rho(\kx(\alpha')_{t_i})=1$.
This proves the claim.

\bigskip
Let us now replace $\alpha$ by a general (in the above sense) $\alpha'\in U$.
Analogously, replace $\gamma=((1+\varepsilon)\alpha-\beta)/\varepsilon$
by the K\"ahler class $\gamma'=((1+\varepsilon)\alpha'-\beta)/\varepsilon$.
In other words, from now on
we are in the following situation: $\alpha\in \kc_X$,
$\gamma=((1+\varepsilon)\alpha-\beta)/\varepsilon\in\kk_X$, and there are
points $t_i\in T(\alpha)\cap S_{\beta_i}$ converging to $0$ such
that $\rho(\kx(\alpha)_{t_i})=1$.

Let $\kx:=\kx(\alpha)\to T(\alpha)$ and denote by $\alpha_t$ a $(1,1)$-class
on $\kx_t$ that spans the orthogonal
complement of $P(\kx_t):=(\sigma_t\cdot\IC+\bar\sigma_t\cdot \IC)_\IR$ in
$F(\alpha)\subset H^2(X,\IR)\cong H^2(\kx_t,\IR)$. We may choose $\alpha_t$
depending continuously on $t$ such that $\alpha_0=\alpha$.
Next, let $\gamma_i:=((1+\varepsilon)\alpha_{t_i}-\beta)/\varepsilon$
which is considered as a class on $X$ or as a class on $\kx_{t_i}$
via the isomorphism $H^2(X,\IR)\cong H^2(\kx_{t_i},\IR)$. Since the union
of the K\"ahler cones $\kk_{\kx_t}$
in $\bigcup_t H^{1,1}(\kx_t)_\IR$ is open, $\gamma_i$
is a K\"ahler class on $\kx_{t_i}$ for $i\gg0$ (see also the arguments in
the proof of Fujiki's Theorem \ref{fujiki}).
On the other hand, for $i\gg0$ the class $\beta_i$ on $\kx_{t_i}$ is of type
$(1,1)$ and $q_X(\beta_i)>0$. Theorem \ref{proj} then asserts that
$\kx_{t_i}$ is projective. Since
$H^{1,1}(\kx_{t_i})_\IQ=\beta_i\cdot\IQ$ and $\beta_i\in\kc_{\kx_{t_i}}$
for $i\gg0$, the class $\beta_i$ is ample on
$\kx_{t_i}$. Thus, $\alpha_{t_i}$ is contained in the segment
$[\gamma_i,\beta_i]$ joining the two K\"ahler classes $\gamma_i$ and $\beta_i$
on $\kx_{t_i}$. Since the K\"ahler cone is convex,
$\alpha_{t_i}$ is a K\"ahler class on $\kx_{t_i}$ for $i\gg0$.

All we need from the preceeding discussion is the following statement:
\begin{itemize}\item[] {\it If $X$ is an irreducible symplectic manifold
and $\alpha\in\kc_X$ is general, then there exists a point $t\in T(\alpha)$
such that $\alpha_t$ is K\"ahler on $\kx_t=\kx(\alpha)_t$.}
\end{itemize}

Let us fix this point $t$. We denote the induced marking
$H^2(\kx_t,\IZ)\cong H^2(X,\IZ)\rpfeil{5}{\varphi}\Gamma$
of $\kx_t$ by $\varphi_t$.
Then $\varphi(F(\alpha))=\varphi(P(X)\oplus\alpha\cdot\IR)=\varphi_t(P(\kx_t)\oplus\alpha_t\cdot\IR)=\varphi_t(F(\alpha_t))$. This way one identifies
$T(\alpha)$ with an open subset of the base of the twistor space associated to
$(\kx_t,\alpha_t)$. (Note that for an arbitrary $(1,1)$-class $\delta$ the
space $T(\delta)$ is only locally defined but if $\delta$ is a K\"ahler class
then $T(\delta)$ means the complete base of the twistor space, 
which is a $\IP^1$, see \ref{deformationhk},\ \ref{twistor2}.)
Hence, there are two families $\kx\to T(\alpha)$  and $\kx'\to T(\alpha)$ over
the same base $T(\alpha)$, where the latter
is (an open subset of) the twistor space
of $\kx_{t}\cong\kx'_t$.
Both deformations are endowed with the natural markings
$\varphi$ and $\varphi'$ such that
${\varphi_t'}^{-1}\circ\varphi_t$ is induced by the isomorphism
$\kx_t\cong\kx'_t$.
Moreover, $\alpha'_s:=({\varphi_s'}^{-1}\circ\varphi_s)(\alpha_s)$ is a 
K\"ahler class on $\kx_s'$ for all $s\in T(\alpha)$. In particular, $\alpha_0'$
is a K\"ahler class on $X':=\kx_0'$. It remains to show that $(X',\varphi'):=
(\kx'_0,\varphi'_0)$ cannot be separated from $(X,\varphi)$.

Consider the maximal open subset $V\subset T(\alpha)$ containing $t$ such that
there exists an isomorphism $\kx|_V\cong\kx'|_V$ extending $\kx_t\cong\kx_t'$.
(Use the Local Torelli Theorem \ref{periodmapdef}
in order to show that $V$ is open.)
Denote by $\overline V$ the closure of $V$ in $T(\alpha)$ and let
$\partial V:=\overline V\setminus V$ be its boundary.
Let $s\in \partial V$. Then there exists an effective cycle $\Gamma=Z+
\sum Y_i\subset \kx_s\times\kx'_s$ as in the proof of
\ref{BRMain}. In particular, $\kx_s\leftarrow Z\to\kx_s'$ is a birational
correspondence, the projections $Y_i\to\kx_s$ and $Y_i\to\kx_s'$ have
positive fibre dimension, and $[\Gamma]_*={\varphi_s'}^{-1}\circ\varphi_s$.
Let us assume that in addition $\kx_s'$ neither
contains non-trivial curves nor effective divisors, e.g.\
$H^{1,1}(\kx_s')_\IZ=0$ (cf.\ \ref{quadraticform}). Then,
$\kx_s\cong Z\cong\kx_s'$ and $[Y_i]_*: H^2(\kx_s,\IZ)\to H^2(\kx_s',\IZ)$
is trivial (see Lemma \ref{lem} below). Since
$\alpha'_s=({\varphi_s'}^{-1}\circ\varphi)(\alpha_s)=[\Gamma]_*(\alpha)=[Z]_*(
\alpha_s)$, this yields that $\alpha_s$ is a K\"ahler class on $\kx_s$.
Next we claim that in fact $\Gamma=Z$. Indeed, if we compute the volume
with respect to the K\"ahler class $p_1^*\alpha_s+p_2^*\alpha'_s$, then
$$\begin{array}{rcl}
vol(\Gamma)&=&vol(Z)+\sum vol(Y_i)\\
&=&\displaystyle{\int_Z (p_1^*\alpha_s+p_2^*\alpha_s')^{2n}}+\sum vol(Y_i)\\
&=&\displaystyle{\int_{\kx_s'}([Z]_*(\alpha_s)+\alpha_s')^{2n}}+\sum vol(Y_i)\\
&=&\displaystyle{\int_{\kx_s'}([\Gamma]_*(\alpha_s)+\alpha_s')^{2n}}+\sum vol(Y_i)\\
&=&\displaystyle{\int_{\kx_s'}(2\alpha_s')^{2n}}+\sum vol(Y_i)\\
\end{array}$$
On the other hand, if $s_i\in V$ converges to $s$, then the volume
of the graph $\Gamma_i$ of the isomorphism $\kx_{s_i}\cong\kx_{s_i}'$ converges
to $vol(\Gamma)$. But $vol(\Gamma_i)=\int_{\Gamma_i}(p_1^*\alpha_{s_i}+p_2^*\alpha'_{s_i})^{2n}=\int_{\kx_{s_i}'}(2\alpha'_{s_i})^{2n}$. Hence,
$\sum vol(Y_i)=0$, i.e.\ $\Gamma=Z$. But this would contradict the
maximality of $V$. Thus, if $s\in\partial V$, then $\kx'_s$
either contains non-trivial curves or divisors.
By \ref{twistor2} the set of points $s\in T(\alpha)$ with this property is 
countable. Hence, $\partial V$ is countable. But then $\partial V$ could not
separate two non-empty open subset $V$ and $T(\alpha)\setminus \overline V$.
Hence, $\overline V=T(\alpha)$, which proves that $(X,\varphi)$ and $(X',\varphi')$ cannot be separated.\qed

\begin{remark}\label{topoo}---
In fact, we can modify $\kx'$ appropriately such that it becomes isomorphic
to $\kx$ over $T(\alpha)\setminus\{0\}$. (However, it will not be
a twistor space any longer.) Indeed, any countable closed subset in
$T(\alpha)$ has an isolated point, but in the neighbourhood
of an isolated point $\ne0$ one can replace $\kx'$ by $\kx$.
\end{remark}

\begin{remark}\label{Weyl}---
The transitivity of the action of the Weyl-group on the set of chambers of
a K3 surface is equivalent to the following statement: If $X$ is a K3 surface
and $\alpha\in\kc_X$ is general in the sense that it is not orthogonal to
any $(-2)$-curve, then there exists a cycle
$\Gamma=\Delta+\sum C_i\times C_i\subset X\times X$,
where $\Delta$ is the diagonal and the $C_i$'s are $(-2)$-curves, such that
$[\Gamma]_*(\alpha)$ is a K\"ahler class. In this light, the above
proposition is a weak generalization of the transitivity of the action
of the Weyl-group. Unfortunately, it seems to be hard to
specify the assumption on $\alpha$ to be `general'. In particular,
one would like to replace it by an open condition. For a K3 surface this 
is granted by the fact that the union of all walls is closed (and hence its
complement is open) in the positive cone. However, as we will see below
(Corollary \ref{Delignehelps}),
a class is `general' in the sense of the proposition
if it is not orthogonal to any integral class. 
\end{remark}

The following lemma is rather elementary. It was used in the previous
proof and will come up again in Sect.\ \ref{kaehlercone}.

\begin{lemma}\label{lem}--- Let $X$ and $X'$ be compact complex manifolds of
dimension $m$ and let $Y\subset X\times X'$ be an $m$-dimensional irreducible
subvariety such that the projection $p':Y\to X'$ has
positive fibre dimension.
If the induced homomorphism $[Y]_*:=p'_*([Y].p^*(\,.\,)):
H^2(X,\IZ)\to H^2(X',\IZ)$ is non-zero, then $p'(Y)\subset X'$ is a divisor.
In this case $[Y]_*(\alpha)=(\int_C\alpha)[p'(Y)]$, where $C={p'}^{-1}(x)$
is the general fibre curve of $p':Y\to p'(Y)$.\qed
\end{lemma}

Before approaching the K\"ahler cone of an arbitrary irreducible symplectic
manifold in Sect.\ \ref{kaehlercone} let us deduce here
some immediate consequences of Proposition
\ref{Weylprop}.

\begin{corollary}\label{Kaehlerconeforperiod}---
Let $X$ be an irreducible symplectic manifold without effective divisors.
Then for general $\alpha\in\kc_X$ there exists
a birational correspondence $Z\subset X\times X'$ between $X$ and another
irreducible symplectic manifold $X'$ with $[Z]_*(\alpha)\in\kk_{X'}$.
If in addition $X$ contains no rational curves, then $\kc_X=\kk_X$
\end{corollary}

\prf If $X$ does not contain any divisor, then by the lemma the maps $[Y_i]:
H^2(X',\IZ)\to H^2(X,\IZ)$ are trivial, where $\Gamma=Z+\sum Y_i$ is as in 
Corollary \ref{corBRMain}. Hence, if $\alpha':=[\Gamma]_*(\alpha)$, which
is a K\"ahler class by Corollary \ref{corBRMain}, then $\alpha=[\Gamma]_*(\alpha')=
[Z]_*(\alpha')$ and, therefore, $[Z]_*(\alpha)=\alpha'\in\kk_{X'}$.
If $X$ does not contain any rational curves, then any birational correspondence
extends to an isomorphism. Hence, $\alpha$ is K\"ahler.\qed

\bigskip
From the last assertion and \ref{twistor2} one easily obtains
also the following 

\begin{corollary}\label{kaehlerconecor}\label{Kaehlerconeforperiod2}---
If $H^{1,1}(X)_\IZ=0$, i.e.\ there exists no non-trivial
line bundle on $X$, then $\kc_X=\kk_X$.
In particular, if $X$ is a general irreducible symplectic manifold
then $\kc_X=\kk_X$.
\end{corollary}

\prf The assumption $H^{1,1}(X)_\IZ=0$ implies
$H^{2n-1,2n-1}(X)_\IZ=0$ by \ref{quadraticform}. Hence,
$X$ contains neither curves nor effective divisors. Thus we can apply
the previous corollary.
For the second assertion recall that the deformations of $X$ that
admit non-trivial line bundles form a countable union of hypersurfaces
in $\Def(X)$ (cf.\ \ref{deformationlb}).\qed

\bigskip
Once the K\"ahler cone of a general irreducible symplectic
manifold is described (Corollary \ref{kaehlerconecor}), 
Proposition \ref{Weylprop} can be sharpened to yield:

\begin{corollary}\label{Delignehelps}
--- Let $(X,\varphi)\in\gM_\Gamma$ be a marked irreducible symplectic manifold.
Assume that $\alpha\in\kc_X$ is not orthogonal to any
$0\ne\beta\in H^{1,1}(X)_\IZ$.
Then there exists a point $(X',\varphi')\in\gM_\Gamma$ which cannot be
separated from $(X,\varphi)$ such that $({\varphi'}^{-1}\circ\varphi)(\alpha)$
is in the K\"ahler cone of $X'$.
\end{corollary}

\prf As in the proof of Proposition \ref{Weylprop} we consider the
`twistor space' $\kx(\alpha)\to T(\alpha)$. By the assumption on $\alpha$,
the general fibre $\kx_t:=\kx(\alpha)_t$ satisfies $H^{1,1}(\kx_t)_\IZ=0$.
By Corollary \ref{kaehlerconecor} this implies $\kk_{\kx_t}=\kc_{\kx_t}$.
If $\alpha_t$ is a $(1,1)$-class on $\kx_t$ that spans the orthogonal
complement of $P(\kx_t):=(\sigma_t\cdot\IC+\bar\sigma_t\cdot \IC)_\IR$ in
$F(\alpha)\subset H^2(X,\IR)\cong H^2(\kx_t,\IR)$ (as in the proof of
\ref{Weylprop}), then the class $\alpha_t$ is K\"ahler
on $\kx_t$ for general $t$ close to $0$. Now proceed as in the final
paragraph of the proof of \ref{Weylprop}\qed.
 
\bigskip
Once a point $(X',\varphi')\in \gM_\Gamma$
non-separated from $(X,\varphi)$ with $({\varphi'}^{-1}\circ\varphi)(\alpha)$
in the K\"ahler cone $\kk_{X'}$ has been  shown to exist, one easily proves

\begin{theorem}\label{Delignesuggestion}---
Let $\kx\to \Def(X)$ be the universal deformation
of an irreducible symplectic manifold $X$ and let $\alpha\in\kc_X$
be a class not orthogonal to any $0\ne\beta\in H^{1,1}(X)_\IZ$.
Then there exists another irreducible symplectic manifold
$X'$ and a smooth proper
family $\kx'\to \Def(X)$ with $X'=\kx'_0$
such that over an open
subset containing the complement of the union of all
hypersurfaces $S_\beta$ with $\beta\in H^2(X,\IZ)$ both
families $\kx$ and $\kx'$ are isomorphic and the induced isomorphism
$H^2(X,\IR)\cong H^2(X',\IR)$ maps $\alpha$ to a K\"ahler class on $X'$.
\end{theorem}

\prf Again, this is just a reformulation of the proof of
\ref{Weylprop} with the extra input \ref{Delignehelps},
which in turn is based on \ref{Weylprop} in its original form.
Fix a marking $\varphi$ of $X$ and take $(X',\varphi')$
as in \ref{Delignehelps}. Consider $\Def(X')$ as an
open subset of $\gM_\Gamma$. Using the period map,
both spaces $\Def(X)$ and $\Def(X')$ can be identified.
Thus, we have two marked families $(\kx,\varphi)\to\Def(X)$
and $(\kx',\varphi')\to\Def(X)$,
both universal for $X$ respectively $X'$.
Moreover, for some point $t\in\Def(X)$ there is an
isomorphism $f:\kx_t\cong\kx'_t$ with $f^*={\varphi'}^{-1}\circ\varphi$.
Let $V$ be the maximal open subset
such that there exists an isomorphism $\kx|_V\cong\kx'|_V$ extending
$f$. As we have seen in the proof of
\ref{Weylprop}, for a point $s$ in
the boundary $\partial V:=\overline V\setminus V$  the group
$H^{1,1}(\kx_s)_\IZ$ is non-trivial.
Hence, $\partial V$ is contained in the union of all hypersurfaces
$S_\beta$ with $\beta\in H^2(X,\IZ)$. Since a countable union
of real codimension two subsets cannot separate two non-empty open subsets,
one of the open sets $V$ or $\Def(X)\setminus \overline V$ must be
empty. Hence, $\overline V=\Def(X)$.\qed

\bigskip
It seems likely that in the above theorem one can achieve that the open subset
is in fact the complement of finitely many hypersurfaces $S_{\beta_i}$
with $\beta_i\in H^{1,1}(X)_\IZ$,
but at the moment I do not know how to prove this.


\section{The Ample Cone}\label{amplecone}

After having established a projectivity criterion for irreducible symplectic
manifolds in Sect.\ \ref{projectivity},
we now strive for an understanding of the ample cone
of a projective irreducible symplectic manifold. Recall the following
well-known results:

\bigskip
\refstepcounter{theorem}\label{Nakai}{\bf \thetheorem} --- Let $X$ be a
projective variety. A line bundle $L$ on $X$ is ample if and only if
$\int_Zc_1^i(L)>0$ for any integral subscheme $Z\subset X$ of dimension
$i\leq\dim(X)$ (cf.\ \cite{HartshorneAmple}).

\refstepcounter{theorem}\label{ampleK3}{\bf \thetheorem} --- Let $X$ be a K3
surface. A line bundle $L$ on $X$ is ample if and only if $c_1(L)$ is
contained in the positive cone and $\int_Cc_1(L)>0$ for all $(-2)$-curves
(i.e.\ irreducible smooth and rational) $C\subset X$ (cf.\ \cite{BPV},
\cite{Periodes}).

\bigskip
The second statement, a special case of the Nakai-Moishezon criterion,
says that on a two-dimensional irreducible
symplectic manifold the ample cone is the integral part
of the positive cone that is positive on all $(-2)$-curves. The main result
of this section (Corollary \ref{ampleconethm}) is formulated in this spirit.
However, the result is extremely weak compared to \ref{ampleK3}
as it only says that an integral class in the positive
cone that cannot be separated from the K\"ahler cone
by any {\it integral} wall is ample.
The problem one has to face in higher dimensions
is that an analogue of $(-2)$-curves is not (yet?) available.

\bigskip
Recall, the following notation: $\kk_X\subset H^{1,1}(X)_\IR$ is the cone
of all K\"ahler classes and the positive cone $\kc_X$ is the component of
$\{\alpha\in H^{1,1}(X)_\IR|q_X(\alpha)>0\}$ that contains $\kk_X$.

We begin with a variant of \ref{Nakai} in the case of projective manifolds
with trivial canonical bundle. The following
proposition, which is a straightforward
consequence of the Basepoint-Free Theorem \cite{CKM}, says
that in order to check whether a line bundle is ample it suffices
to test it on subvarieties which are either $X$ itself or of dimension one.

\begin{proposition}---
Let $X$ be a projective manifold with $K_X\cong\ko_X$. Then a line bundle
$L$ on $X$ is ample if and only if
$\int c_1(L)^{\dim(X)}>0$ and $\int_Cc_1(L)>0$ for all
curves $C\subset X$.
\end{proposition}

\prf Obviously, if $L$ is ample, then both inequalities hold.
If $L$ is a line bundle satisfying both inequalities,
then $L$ is nef and $L\otimes K\dual_X\cong L$ is nef and big.
Then, the Basepoint-Free Theorem shows that $L^m$ is globally generated
for $m\gg0$.
But any globally generated line bundle that is positive on all curves is ample
(cf.\ \cite{HartshorneAmple}).\qed

\bigskip
Applied to irreducible symplectic manifolds this yields:

\begin{corollary}\label{NakaionHK}--- 
Let $X$ be a projective irreducible symplectic manifold
of dimension $2n$ and let
$L$ be a line bundle on $X$. Then
\begin{itemize}
\item[$\bullet$] $L$ is ample if and only if $c_1(L)\in\kc_X$ and $\int_Cc_1(L)>0$
for all curves $C\subset X$.
\item[$\bullet$] $L$ is in the closure of the ample cone if and only if
$c_1(L)\in\overline\kc_X$ and $L$ is nef.
\end{itemize}
\end{corollary}

\prf Observe that for the first assertion we do not need the projectivity
of $X$, for $c_1(L)\in\kc_X$ implies
that $X$ is projective (Theorem \ref{proj}). Both assertions follow
directly from the proposition: If $L$ is ample, then
certainly $\int_Cc_1(L)>0$, $q_X(c_1(L))>0$, and
$q_X(c_1(L),\alpha)>0$ for any
K\"ahler class $\alpha$ (cf.\ \ref{posconedef})
and, therefore, $c_1(L)\in\kc_X$.
Conversely, if $c_1(L)\in\kc_X$ then $\int_Xc_1^{2n}(L)>0$ by
\ref{Todd}.
For the second assertion replace $L$ by $L^m\otimes M$ with $M$ an
ample line bundle and $m\gg0$.\qed

\bigskip
Using the isomorphism $cL_{\sigma\bar\sigma}^{n-1}:H^{1,1}(X)_\IQ\cong
H^{2n-1,2n-1}(X)_\IQ$ (cf.\ \ref{quadraticform}) Corollary \ref{NakaionHK}
allows us to formulate also the following

\begin{corollary}\label{ampleconethm}--- 
Let $X$ be an irreducible symplectic manifold and let $L$ be a 
line bundle on $X$. Then $L$ is ample if and only if
$c_1(L)$ satisfies
\begin{itemize}
\item[{\it i)}] $q_X(c_1(L),\,.\,)$ is positive on $\kk_X$.
\item[{\it ii)}] If $M\in \Pic(X)$ such that $q_X(c_1(M),\,.\,)$ is positive on
$\kk_X$, then $q_X(c_1(M),c_1(L))>0$.
\end{itemize}
\end{corollary}

\prf If $L$ is ample, then $c_1(L)\in\kk_X$ and, therefore, {\it i)} and
{\it ii)} follow. If $L$ satisfies {\it i)} and {\it ii)}, then $c_1(L)\in\kc_X$. Hence, $X$ is projective (Theorem \ref{proj}).
Thus it remains to verify that $L$ is positive on all curves
(\ref{NakaionHK}).
If $C$ is a curve then $[C]\in H^{2n-1,2n-1}(X)_\IQ$ and
therefore there exists an $\alpha\in H^{1,1}(X)_\IQ$ with
$cL_{\sigma\bar\sigma}^{n-1}(\alpha)=[C]$. Since any K\"ahler
class $[\omega]$ is positive on $C$, one has $q_X(\alpha,\omega)>0$.
Hence, $q_X(\alpha,c_1(L))>0$ by assumption, for $\alpha$
is a rational class.\qed

\bigskip
The statement seems rather weak, and indeed, it cannot be considered as
a true generalization of the Nakai-Moishezon criterion
for K3 surfaces. However, it is non-trivial. E.g.,
if $X$ is a projective irreducible symplectic manifold such that its Picard
group is spanned by two line bundles $L_1$ and $L_2$ with $L_1$ ample
and $L_2$ non-ample, then there exists a ($\IZ$-)linear combination
of $L_1$ and $L_2$ that is negative (with respect to $q_X$)
on $L_2$ and positive on the whole K\"ahler (and not only on the ample!) cone. 
Of course, if $h^{1,1}(X)=\rho(X)$, then the assertion of the corollary
is void.

Note that the description of the ample cone given in this section
does not make use of Proposition \ref{Weylprop}, but 
I expect that Corollary \ref{ampleconethm} together with
Proposition \ref{Weylprop} implies that the ample cone
is finitely polyhedral. But at the moment I do not
know how to prove this.

\section{The K\"ahler Cone}\label{kaehlercone}

Here we slightly generalize the result of the previous section
and give a description of the K\"ahler cone. Again,
the result is much weaker than the known ones for K3 surfaces and
says that a class in the positive cone that cannot be separated
by any integral wall from the K\"ahler cone is K\"ahler itself.
It might be noteworthy that the description of the ample cone
of a K3 surface is a rather easy consequence of the Nakai-Moishezon criterion,
whereas the description of the K\"ahler cone of
a K3 surface relies on the Global Torelli Theorem,
which is in no form available in higher dimensions.
 
The result can most powerfully be applied to the case where there
are no or only few integral classes. In particular it generalizes
Corollary \ref{kaehlerconecor} to the case that $X$ is projective and
the Picard number is one (cf.\ \ref{nameless}).

\begin{theorem}\label{kaehlerconethm}---
Let $X$ be an irreducible symplectic manifold. Then a class
$\alpha$ is contained in the closure of the K\"ahler cone $\overline{\kk}_X$
if and only if
\begin{itemize}
\item[{\it i)}] $\alpha\in\overline{\kc}_X$.
\item[{\it ii)}] If $q_X(c_1(M),\,.\,)$ is non-negative on the K\"ahler cone
$\kk_X$ for a line bundle
$M\in\Pic(X)$, then $q_X(c_1(M),\alpha)\geq0$.
\end{itemize}
\end{theorem}

\prf The principal idea of the proof is modelled on Beauville's
expos\'e in \cite{Periodes}. Arguments using the Weyl-group
are replaced by Proposition \ref{Weylprop}.

It is obvious that {\it i)} and {\it ii)} are necessary
for $\alpha\in\overline{\kk}_X$. For the converse,
assume that $\alpha$ satisfies
both conditions. We first consider $\alpha+\varepsilon\cdot\gamma$,
where $0<\varepsilon\ll1$ and $\gamma$ is a general K\"ahler class.
This class is contained in $\kc_X$ and satisfies the strong inequality
in {\it ii)} for all line bundles $M$ with $q_X(c_1(M),\,.\,)$ positive
on $\kk_X$. If we can show that for general $\gamma$ the class
$\alpha+\varepsilon\gamma$ is a K\"ahler class, then
$\alpha\in\overline{\kk}_X$.
Certainly, $\alpha+\varepsilon\gamma$ is general (in the sense of
Proposition \ref{Weylprop}) if $\gamma$ is a general K\"ahler class.
Thus, we only have to deal with the following situation: $\alpha$
is a general class in $\kc_X$ such that $q_X(c_1(M),\alpha)>0$ whenever
$q_X(c_1(M),\,.\,)$ is positive on $\kk_X$, where $M$ is an arbitrary line
bundle on $X$. In this situation we show that $\alpha$ is a K\"ahler class.

Indeed, by Proposition \ref{Weylprop} and Corollary
\ref{corBRMain} there exists another irreducible
symplectic manifold $X'$ together with an effective cycle
$\Gamma=Z+\sum Y_i\subset X\times X'$ satisfying
{\it i)-iv)} of Corollary \ref{corBRMain}.
In particular, there exists a K\"ahler class $\alpha'$ on $X'$
with $[\Gamma]_*(\alpha)=\alpha'$.
The latter condition is equivalent to $[\Gamma]_*(\alpha')=\alpha$.
Since $[\Gamma]_*$ respects the quadratic form,
we have $q_X(\alpha)=q_{X'}(\alpha')$.
By \ref{hodgeunderbir}
any birational correspondence $Z$ respects the quadratic form
as well. Hence, $q_{X'}(\alpha')=q_{X}([Z]_*(\alpha'))$.
This yields
$$\begin{array}{ccl}
0&=&q_X(\alpha)-q_X([Z]_*(\alpha'))\\
&=&q_X(\alpha+[Z]_*(\alpha'),\alpha-[Z]_*(\alpha'))\\
&=&q_X(\alpha+[Z]_*(\alpha'),\sum[Y_i]_*(\alpha')).\\
\end{array}$$
By Lemma \ref{lem} only those components $Y_i$ contribute
for which $V_i:=p(Y_i)$ is a divisor. Moreover, in this case
$[Y_i]_*(\alpha')=(\int_{C_i}\alpha')[V_i]$, where $C_i$ is the generic fibre
of $p:Y_i\to V_i$. Since $\alpha'$ is K\"ahler, $\int_{C_i}\alpha'>0$.
By \ref{Todd} a K\"ahler class
is positive (with respect to $q$) on any effective divisor. Hence,
any $q_X((\int_{C_i}\alpha')[V_i],\,.\,)$ is positive on $\kk_X$.
Using the assumption on $\alpha$, this implies
$q_X(\sum(\int_{C_i}\alpha')[V_i],\alpha)>0$.
Since $\alpha'$ is K\"ahler, one also has
$q_X([Z]_*(\alpha'),[V_i])=q_{X'}(\alpha',[Z]_*[V_i])>0$
(cf.\ \ref{Todd}).
Altogether, this yields that if $\sum[Y_i]_*(\alpha')\ne0$ then
$0=q_X(\alpha+[Z]_*(\alpha'),\sum(\int_{C_i}\alpha')[V_i])>0$.
Hence, $\sum[Y_i]_*(\alpha')=0$ and, therefore,
$\alpha=[\Gamma]_*(\alpha')=[Z]_*(\alpha')$.

If $Z$ does not define an
isomorphism $X\cong X'$ then there exists a (rational) curve $C\subset X$
such that $\int_C\alpha=\int_C[Z]_*(\alpha')<0$, e.g.\
take a curve in the fibre
of $Z\to X'$. On the other hand, there exists a rational class
$\beta\in H^{1,1}(X)_\IQ$ such that $cL_{\sigma\bar\sigma}^{n-1}(\beta)=
[C]$ (cf.\ \ref{quadraticform}). Since any K\"ahler class is positive on
$C$, the linear form $q_X(\beta,\,.\,)$ is positive on
$\kk_X$. By the assumption this yields $\int_C\alpha>0$.
Contradiction. Thus $X\cong Z\cong X'$ and, therefore, $\alpha\in\kk_X$.\qed

\bigskip
The following is another instance where the K\"ahler cone can completely
be described in terms of the period.

\begin{corollary}\label{nameless}---
Let $X$ be an irreducible symplectic manifold and assume that $\Pic(X)$ is
spanned by a line bundle $L$ such that $q_X(c_1(L))\geq0$.
Then $\kk_X=\kc_X$, i.e.\ any class $\alpha\in\kc_X$ is K\"ahler.
\end{corollary}

\prf Of course, it suffices to prove that $\overline{\kc}_X=\overline{\kk}_X$.
Thus we can apply the theorem.
Without loss of generality we can assume that
$q_X(c_1(L),\,.\,)$ is non-negative
on $\overline{\kc}_X$ (Hodge Index).
Therefore, $\overline{\kc}_X\subset\overline{\kk}_X$.\qed


\section{Surjectivity of the Period Map}\label{periodmap}

Recall, that for a lattice $\Gamma$ of index $(3,b-3)$
we defined the period domain $Q$ as
the set $\{x\in\IP(\Gamma_\IC)|q_\Gamma(x)=0,~q_\Gamma(x+\bar x)>0\}$,
which is an open set of a smooth quadric. Also recall,
that the period map $\kp:\gM_\Gamma\to \IP(\Gamma_\IC)$ takes values
in $Q$. (cf.\ \ref{moduli}).
In this section we present a proof of the following

\begin{theorem}\label{surjper}---
Let $\gM_\Gamma^o$ be a non-empty connected component of the
moduli space $\gM_\Gamma$ of marked irreducible symplectic manifolds. Then the
period map
$$\kp:\gM_\Gamma^o\to Q$$
is surjective.
\end{theorem}

The proof is, once again, modelled on Beauville's presentation in
\cite{Periodes}. One proceeds in two steps. The first part of the proof
consists of showing that all points of $Q$ are equivalent with respect to
a certain equivalence relation defined below. This part is a word-by-word
copy of the known arguments.
The second part, where it is shown
that the image of the period map is invariant under the equivalence
relation, deviates from the standard proofs even for K3 surfaces.
The description of the K\"ahler cone (Corollary \ref{Kaehlerconeforperiod})
turns out to be crucial
for this part.

Let us first recall the following lemma (cf.\ \cite{Periodes}):

\begin{lemma}--- The map sending a point $x\in Q\subset\IP(\Gamma_\IC)$ to
$P(x):=(x\cdot\IC\oplus\bar x\cdot\IC)\cap\Gamma_\IR$
defines a natural isomorphism between $Q$ and the Grassmannian $\tilde Q$
of positive oriented planes in $\Gamma_\IR$.\qed
\end{lemma}

This lemma enables one to prove that the period domain $Q$ is connected.
See \cite{Periodes} for the complete argument.

\begin{definition}---
Two points $x,y\in Q$ are called equivalent if there exists a sequence
$x=x_1,x_2,\ldots,x_k=y\in Q$ such that the subspaces
$\langle P(x_i),P(x_{i+1})\rangle\subset\Gamma_\IR$
are of dimension three and such that
$q_\Gamma|_{\langle P(x_i),P(x_{i+1})\rangle}$ is positive definite.
(A subspace of dimension three with this property is called a positive
$3$-space.)
\end{definition}

One can easily show that the set of points equivalent to a fixed $x\in Q$
is open. According to \cite{Periodes} this
together with the connectivity of
the period domain $Q$ is enough to prove
(again, for the complete argument see \cite{Periodes}):

\begin{lemma}---
Any two points in $Q$ are equivalent.\qed
\end{lemma} 

That these results are valid in this generality, and hence
applicable to higher dimensional irreducible symplectic manifolds,
was also noticed by Verbitsky
\cite{Verbitsky}.

This lemma together with the following one immediately proves
Theorem \ref{surjper}.

\begin{lemma}---
If $x,y\in Q$ such that $\langle P(x),P(y)\rangle$ is a positive $3$-space,
then $x\in\kp(\gM_\Gamma^o)$ if and only if $y\in\kp(\gM_\Gamma^o)$.
\end{lemma}

\prf Assume $x=\kp((X,\varphi))$, i.e.\ $P(x)=\varphi(P(X))$,
with $(X,\varphi)\in\gM_\Gamma^o$.

I claim that one can deform $X$ slightly such that
\begin{itemize}
\item[{\it i)}] $\langle P(x),P(y)\rangle$ is still a positive $3$-space.
\item[{\it ii)}] $\rho(X)=0$.
\end{itemize}

This is proved as follows: Identify $\Def(X)$ with a small
open neighbourhood of $x$ in $Q$ via the period map (Local Torelli
Theorem). Then consider the countable union of hypersurfaces
$T=\bigcup S_\alpha\subset \Def(X)$ where
$\alpha\in H^2(X,\IZ)$ (see
\ref{hypersurfaces}).
To achieve {\it i)} and {\it ii)} it suffices to find
$t\in \Def(X)\setminus T$ such that $\langle \varphi(P(\kx_t)),P(y)\rangle$
is a positive $3$-space. Of course, the positivity is harmless as long
as $t$ is close to $0$ and $\dim\langle \varphi(P(\kx_t)),P(y)\rangle=3$.

First, consider those $t$ for which
$\langle \varphi(P(\kx_t)),P(y)\rangle=\langle P(x),P(y)\rangle$. They are
parametrized by (an open subset of) of the non-degenerate
quadric $Q\cap \langle P(x),P(y)\rangle_\IC$ (for a similar argument see
\ref{twistor2}). Moving $t$ slightly in $Q\cap\langle P(x),P(y)\rangle_\IC$
we can assume that the orthogonal complement $k_0\cdot\IR\subset\langle P(x),
P(y)\rangle$ of $P(x)$ is not contained in $P(y)$. Now fix a basis
$P(y)=\langle a,b\rangle$. Then to any $k$ in a neighbourhood 
of $k_0$ in $\Gamma_\IR$ we associate
$P_k:=\langle q(k,k)a-q(a,k)k,q(k,k)b-q(b,k)k\rangle$ -- the orthogonal
complement of $k$ in $\langle k,P(y)\rangle$. Then $\langle P_k,P(y)\rangle$
is a positive $3$-space.

Let $T'\subset \Def(X)$ be the subset of those $t\in\Def(X)$ such that
$\varphi(P(\kx_t))=P_k$ for some $k$ in a neighbourhood of $k_0\in\Gamma_\IR$.
If $T'\not\subset T$ we are done. If $T'\subset T$, then there exists
an $\alpha$ such that $T'\subset S_\alpha$ (note that $T'$ is locally
irreducible).
Since $T'$ certainly contains those points
for which $\langle\varphi(P(\kx_t)),P(y)\rangle=\langle P(x),P(y)\rangle$,
we in particular have (locally!) $Q\cap\langle P(x),P(y)\rangle_\IC\subset
T'\subset S_\alpha$. Since $S_\alpha$ is the hyperplane section
defined by $q(\alpha,\,.\,)$ and $Q\cap \langle P(x),P(y)\rangle_\IC$ is a 
non-degenerate quadric, this shows that $q(\alpha,\,.\,)$ vanishes
on $\langle P(x),P(y)\rangle$. Moreover, since $T'\subset S_\alpha$,
the spaces $P_k$ are all orthogonal to $\alpha$.
Hence, $q(q(k,k)a-q(a,k)k,\alpha)=0$ and  $q(q(k,k)b-q(b,k)k,\alpha)=0$
for all $k$ in a neighbourhood of $k_0$.
Also $q(a,\alpha)=q(b,\alpha)=0$ and, therefore, $q(a,k)q(k,\alpha)
=q(b,k)q(k,\alpha)=0$ for all $k$ in a neighbourhood of $k_0$. Contradiction.

Corollary \ref{Kaehlerconeforperiod} then shows that for
$X$ satisfying {\it i)} and {\it ii)} above one has $\kc_X=\kk_X$. Thus,
if $\beta$ spans the orthogonal complement of $P(x)$ in
$\langle P(x),P(y)\rangle$, then,
after replacing $\beta$ by $-\beta$ if necessary,
$\alpha:=\varphi^{-1}(\beta)\in\kc_X$ is a
K\"ahler class. On the other hand, $y\in\IP(\varphi(F(\alpha)_\IC))\subset
\IP(\Gamma_\IC)$. Since there exists the twistor space $\kx(\alpha)\to
T(\alpha)=\IP(\varphi(F(\alpha)_\IC))$ and $T(\alpha)\subset\gM_\Gamma^o$,
this suffices to conclude that $y\in \kp(\gM_\Gamma^o)$.\qed

\bigskip
Note that in the original proof for K3 surfaces the class $\beta$ is chosen
to be rational and not orthogonal to any $(-2)$-curve. Therefore, modulo
the action of the Weyl-group, it is positive on all $(-2)$-curves and, hence,
ample by the Nakai-Moishezon criterion. In the approach above the
generalization of the Nakai-Moishezon criterion as proved in Sect.\
\ref{amplecone} does not suffice to prove the surjectivity along this line.
Also note, that the order of the arguments is different compared to the
original proof for K3 surfaces. For K3 surfaces one first proves the
surjectivity of the period map using the Nakai-Moishezon criterion
and then applies it to derive a description of
the K\"ahler cone. Whereas here, we first described the K\"ahler cone and
then applied the result, which is, however, far less explicit than
the known one for K3
surfaces, to prove the surjectivity of the period map.

\section{Automorphisms}\label{auto}

Not much is known about the automorphism group of an irreducible
symplectic manifold of dimension greater than two.
In this section we collect some results related to
this question and add some remarks.

We introduce the following notations: Let $X$ be an irreducible
symplectic manifold. Then
\begin{itemize}
\item[--] $Aut(X)$ is the group of holomorphic automorphisms of $X$.
\item[--] $Birat(X)$ is the group of birational
automorphisms, i.e.\ of birational maps $X- - \to X$.
\item[--] $A(X)\subset Aut(H^2(X,\IZ))$ is the group of automorphisms of
$H^2(X,\IZ)$ which are com\-pa\-tible with the Hodge structure and the quadratic
form $q_X$ and map the K\"ahler cone to the K\"ahler cone.
\item[--] $B(X)\subset Aut(H^2(X,\IZ))$ is the group
of automorphisms of $H^2(X,\IZ)$ which are
compa\-tible with the Hodge structure and
the quadratic form $q_X$.
\item[--] $a:Aut(X)\to Aut(H^2(X,\IZ))$ maps an automorphism
$f$ to $f^*$.
\item[--] $b:Birat(X)\to Aut(H^2(X,\IZ))$ maps a birational map
$f$ to $f^*$.
\end{itemize}
\begin{proposition}---
{\it i)} $\im(a)\subset A(X)$, {\it ii)} $\im(b)\subset B(X)$,
{\it iii)} $b^{-1}(A(X))=Aut(X)\subset Birat(X)$, {\it iv)}
$\ker(a)=\ker(b)$, and {\it v)} $\ker(a)$ is finite. 
\end{proposition}

\prf {\it i)} is trivial, {\it ii)} follows from Lemma \ref{hodgeunderbir}.
In order to prove {\it iii)} one has to show that
any birational map which maps a K\"ahler class to a K\"ahler class
can be extended to an automorphism.  This was proved in \cite{Fujiki3}.
{\it iv)} is an easy consequence of {\it iii)}. To prove {\it v)}
one evokes two standard facts: Firstly, the group of isometries of a compact
Riemannian manifold is compact and, secondly, the Calabi-Yau metric
with respect to a fixed K\"ahler class is unique. Hence, an automorphism
acting trivially on the second cohomology leaves invariant the K\"ahler
class and, hence, the Calabi-Yau metric, i.e.\ it is an isometry.
Since the group $Aut(X)$ is discrete, this suffices to conclude
that $\ker(a)$ is finite. \qed

\bigskip
The natural inclusion  $Aut(X)\subset Birat(X)$ is in general
proper. Indeed, Beauville constructed in \cite{Beauville2} an example
of a birational automorphism of $\Hilb^n(S)$, where $S$ is a
special K3 surface, which does not extend to an automorphism (cf.\
\ref{bevdeb}).

However, one has

\begin{proposition}---
If $X$ is general, i.e.\ together with a marking $\varphi$
it is a general point in $\gM_\Gamma$, then $Aut(X)=Birat(X)$.
\end{proposition}

\prf This stems from the fact that the general irreducible symplectic manifold
$X$ does not contain any (rational) curves (cf.\ \ref{twistor2}).
But a birational
automorphism of a variety without rational curves extends to a
holomorphic automorphism.\qed

\bigskip
For K3 surfaces one certainly has $Aut(X)=Birat(X)$. Moreover, the Global
Torelli Theorem for K3 surfaces in particular asserts that $a:Aut(X)\to A(X)$
is an isomorphism. Using this fact one also shows that
$a$ is injective for the Hilbert scheme of a K3 surface
\cite{Beauville2}. However, due to an example
of Beauville \cite{Beauville2}, we know
that $a$ is not injective in general (This time the counterexample is provided
by a generalized Kummer variety). Since the graph of an automorphism in the
kernel of $a$ deforms (at least infinitesimally) in all directions with $X$,
we shall not even expect that $a$ is injective for general $X$.

\bigskip
{\bf Questions} --- {\it i)} Is $a$ surjective? {\it ii)}
Is $A(X)=\{1\}$ for general $X$?

\bigskip
The fact that $a$ is not injective has the unpleasant consequence
that the moduli space $\gM_\Gamma$ of marked irreducible symplectic manifolds
is not fine (contrary to the K3 surface case). For
a discussion of this and other questions related to the moduli
space see the next section.


\section{Further Remarks}\label{remarks}

As mentioned in the introduction, the two problems I consider
the most important in the theory are:
\begin{itemize}
\item[--] Is there a Global Torelli Theorem for irreducible symplectic
manifolds in higher dimensions?
\item[--] What are the possible deformation (or diffeomorphism) types of
irreducible symplectic manifolds?
\end{itemize}

The second question alludes to the rather easy fact that any two K3 surfaces
are deformation equivalent and, hence, diffeomorphic. From dimension
four on we know exactly two different deformation (diffeomorphism)
types of irreducible symplectic manifolds (cf.\ Sect.\ \ref{biratman}).
But there is no obvious reason why there should be no other possibilities.

Due to Debarre's counterexample \cite{Debarre} (cf.\ \ref{bevdeb}) we know
that the the Global Torelli Theorem as formulated for K3 surfaces
fails in higher dimensions.
However, there seems to be no counterexample known to the following
speculation which was first formulated by Mukai \cite{Mu2}.

\bigskip
{\bf Speculation}\refstepcounter{theorem}\label{Spec} {\bf \thetheorem}
--- {\it Two irreducible symplectic manifolds with
isomorphic periods are birational.}

\bigskip
The known proofs of the Global Torelli Theorem for
K3 surfaces break down in higher dimensions.
This is mainly due to a missing analogue of Kummer surfaces.
Kummer surface are dense in the moduli spaces of K3 surfaces and 
whether a K3 surface is a Kummer surface
can easily be read off its period. (Quartic hypersurfaces would
be another such distinguished class of K3 surfaces \cite{Friedman}.)
In higher dimensions we neither have a good class of manifolds
which are dense in $\gM_\Gamma$ nor do we know a `typical' class
of manifolds which could be recognized by its period.
In this light, it would be interesting to find an answer to the
following

\bigskip
{\bf Question} ---
Let $X$ be an irreducible symplectic manifold with the period
of a Hilbert scheme $\Hilb^n(S)$ of a K3 surface $S$. Is
$X$ birational to $\Hilb^n(S)$?

\bigskip
Let us conclude with a few remarks on the moduli space $\gM_\Gamma$
of marked manifolds.

By definition $\gM_\Gamma$ is the set $\{(X,\varphi)\}/\sim$.
Here $X$ is an irreducible symplectic manifold and $\varphi:H^2(X,\IZ)\cong\Gamma$ is an isomorphism of lattices. Two marked manifolds $(X,\varphi)$
and $(X',\varphi')$ are equivalent if there exists an isomorphism
$g:X\cong X'$ with $g^*=\pm({\varphi}^{-1}\circ\varphi')$.
That this set can be given the structure of a smooth, although non-separated,
manifold is a consequence of the unobstructedness of symplectic manifolds
and the Local Torelli Theorem.
The period map $\kp:\gM_\Gamma\to Q\subset\IP(\Gamma_\IC)$
exhibits $\gM_\Gamma$ as a space \'etale over of the period domain $Q$.

The above speculation together with 
\ref{BRMain} and \ref{surjper} is equivalent to

\bigskip
{\bf Speculation}\label{Spec'} {\bf \thetheorem'}
--- {\it If $\gM_\Gamma\ne\emptyset$
then the fibre over a general point $x\in Q$ is
exactly one point.}

\bigskip
Note also that once the first speculation is answered positively it can be used
to generalize Theorem \ref{birat} to the effect that any two birational
irreducible symplectic manifolds (projective or not) correspond
to non-separated points in the moduli space.

It is tempting to define another moduli space which does not distinguish
between birational manifolds. Let $\gN_\Gamma$ be the set
$\{(X,\varphi)\}/\approx$, where the $(X,\varphi)$ are as above and
$(X,\varphi)\approx(X',\varphi')$ if and only if there exists a birational map
$g:X- - \to X'$ with  $g^*=\pm({\varphi}^{-1}\circ\varphi')$.
Of course, there is a natural map $\gM_\Gamma\to\gN_\Gamma$.
Theorem \ref{birat} proves that two points in the same fibre
of this map,
at least if they are projective, are contained in the
same component of $\gM_\Gamma$. In fact, they are non-separated there. Thus,
the general form of Theorem \ref{birat}, i.e.\ without the projectivity
assumption, would prove that the number of components of $\gM_\Gamma$ and
$\gN_\Gamma$ is the same.
Can $\gN_\Gamma$ be endowed with the structure of a manifold?
If a birational map can always be extended to birational maps
between all nearby fibres in the Kuranishi family, then
local patching would supply $\gN_\Gamma$ with the structure of
a manifold. Of course, if the Global Torelli Theorem in the formulation
of the first speculation above holds true, then this would be trivial.
Also note that the period map $\kp:\gM_\Gamma\to Q$ naturally factorizes
through a map $\gN_\Gamma\to Q$.

\end{document}